\documentclass[12pt]{article}%
\usepackage{amssymb}
\usepackage{amsmath}
\usepackage{geometry}
\usepackage[onehalfspacing]{setspace}
\usepackage[dvips]{graphicx}%
\usepackage{amsfonts}
\providecommand{\U}[1]{\protect\rule{.1in}{.1in}}
\newtheorem{theorem}{Theorem}

\newtheorem{conjecture}[theorem]{Conjecture}
\newtheorem{corollary}[theorem]{Corollary}

\newtheorem{lemma}[theorem]{Lemma}

\newtheorem{proposition}[theorem]{Proposition}
\newtheorem{remark}[theorem]{Remark}

\begin{document}

\title{Alternative asymptotics for cointegration tests in large VARs.}
\author{Alexei Onatski and Chen Wang\\Faculty of Economics, University of Cambridge}
\maketitle

\begin{abstract}
Johansen's (1988, 1991) likelihood ratio test for cointegration rank of a
Gaussian VAR depends only on the squared sample canonical correlations between
current changes and past levels of a simple transformation of the data. We
study the asymptotic behavior of the empirical distribution of those squared
canonical correlations when the number of observations and the dimensionality
of the VAR diverge to infinity simultaneously and proportionally. We find that
the distribution almost surely weakly converges to the so-called
\textit{Wachter distribution}. This finding provides a theoretical explanation
for the observed tendency of Johansen's test to find \textquotedblleft
spurious cointegration\textquotedblright. It also sheds light on the workings
and limitations of the Bartlett correction approach to the over-rejection
problem. We propose a simple graphical device, similar to the scree plot, for
a preliminary assessment of cointegration in high-dimensional VARs.

\end{abstract}

\section{Introduction}

Johansen's (1988, 1991) likelihood ratio (LR) test for cointegration rank is a
very popular econometric technique. However, it is rarely applied to systems
of more than three or four variables. On the other hand, there exist many
applications involving much larger systems. For example, Davis (2003)
discusses a possibility of applying the test to the data on seven aggregated
and individual commodity prices to test Lewbel's (1996) generalization of the
Hicks-Leontief composite commodity theorem. In a recent study of exchange rate
predictability, Engel, Mark, and West (2015) contemplate a possibility of
determining the cointegration rank of a system of seventeen OECD exchange
rates. Banerjee, Marcellino, and Osbat (2004) emphasize the importance of
testing for no cross-sectional cointegration in panel cointegration analysis
(see Breitung and Pesaran (2008) and Choi (2015)), and the cross-sectional
dimension of modern macroeconomic panels can easily be as large as forty.

The main reason why the LR test is rarely used in the analysis of relatively
large systems is its poor finite sample performance. Even for small systems,
the test based on the asymptotic critical values does not perform well (see
Johansen (2002)). For large systems, the size distortions become overwhelming,
leading to severe over-rejection of the null in favour of too much
cointegration\ as shown in many simulation studies, including Ho and Sorensen
(1996) and Gonzalo and Pitarakis (1995, 1999).

In this paper, we study the asymptotic behavior of the sample canonical
correlations that the LR statistic is based on, when the number of
observations and the system's dimensionality go to infinity simultaneously and
proportionally. We show that the empirical distribution of the squared sample
canonical correlations almost surely converges to the so-called
\textit{Wachter distribution} which also arises, albeit with different
parameters, as the limit of the empirical distribution of the squared sample
canonical correlations between two independent high-dimensional white noises
(see Wachter (1980)). Our analytical findings explain the observed
over-rejection of the null hypothesis by the LR test, shed new light on the
workings and limitations of the Bartlett-type correction approach to the
problem (see Johansen (2002)), and lead us to propose a very simple graphical
device, similar to the scree plot, for a preliminary analysis of the validity
of cointegration hypotheses in large vector autoregressions.

The basic framework for our analysis is standard. Consider a $p$-dimensional
VAR in the error correction form
\begin{equation}
\Delta X_{t}=\Pi X_{t-1}+%
{\displaystyle\sum_{i=1}^{k-1}}
\Gamma_{i}\Delta X_{t-i}+\Phi D_{t}+\varepsilon_{t},\label{model}%
\end{equation}
where $D_{t}$ and $\varepsilon_{t}$ are vectors of deterministic terms and
zero-mean Gaussian errors with unconstrained covariance matrix, respectively.
The LR statistic for the test of the null hypothesis of no more than $r$
cointegrating relationships between the $p$ elements of $X_{t}$ against the
alternative of more than $r$ such relationships is given by
\begin{equation}
LR_{r,p,T}=-T%
{\displaystyle\sum_{i=r+1}^{p}}
\log\left(  1-\lambda_{i}\right)  \text{,}\label{LR statistic}%
\end{equation}
where $T$ is the sample size, and $\lambda_{1}\geq...\geq\lambda_{p}$ are the
squared sample canonical correlation coefficients between residuals in the
regressions of $\Delta X_{t}$ and $X_{t-1}$ on the lagged differences $\Delta
X_{t-i},$ $i=1,...,k-1,$ and the deterministic terms.

In the absence of the lagged differences and deterministic terms, the
$\lambda$'s are the eigenvalues of $S_{01}S_{11}^{-1}S_{01}^{\prime}%
S_{00}^{-1},$ where $S_{00}$ and $S_{11}$ are the sample covariance matrices
of $\Delta X_{t}$ and $X_{t-1},$ respectively, while $S_{01}$ is the cross
sample covariance matrix. More substantively, $\lambda_{1}$ is the largest
possible squared sample correlation coefficient between arbitrary linear
combinations of the entries of $\Delta X_{t}$ and the entries of $X_{t-1},$
$\lambda_{2}$ is the largest squared correlation among linear combinations
restricted to be orthogonal to those yielding $\lambda_{1}$, and so on (see
Muirhead (1982), ch. 11).

Johansen (1991) shows that the asymptotic distribution of $LR_{r,p,T}$ under
the asymptotic regime where $T\rightarrow\infty$ while $p$ remains fixed, can
be expressed in terms of the eigenvalues of a matrix whose entries are
explicit functions of a $p-r$-dimensional Brownian motion. Unfortunately, for
relatively large $p$, this asymptotics does not produce good finite sample
approximations, as evidenced by the over-rejection phenomenon mentioned above.
Therefore, in this paper, we consider a \textit{simultaneous} asymptotic
regime $p,T\rightarrow_{c}\infty$ where both $p$ and $T$ diverge to infinity
so that
\begin{equation}
p/T\rightarrow c\in\left(  0,1\right]  ,\label{asymptotic regime}%
\end{equation}
while $p$ remains no larger than $T$. Our Monte Carlo analysis shows that the
corresponding asymptotic approximations are relatively accurate even for such
small sample sizes as $p=10$ and $T=20$.

The basic specification for the data generating process (\ref{model}) that we
consider has $k=1$. In the next section, we discuss extensions to more general
VARs with low-rank $\Gamma_{i}$ matrices and additional common factor terms.
We also explain there that our main results hold independently from whether a
deterministic vector $D_{t}$ with fixed or slowly-growing dimension is present
or absent from the VAR.

Our study focuses on the behavior of the empirical distribution function
(d.f.) of the squared sample canonical correlations,%
\begin{equation}
F_{p,T}\left(  \lambda\right)  =\frac{1}{p}%
{\displaystyle\sum_{i=1}^{p}}
\mathbf{1}\left\{  \lambda_{i}\leq\lambda\right\}  ,\label{empirical d.f.}%
\end{equation}
where $\mathbf{1}\left\{  \cdot\right\}  $ denotes the indicator function. We
find that, under the null of $r$ cointegrating relationships, as
$p,T\rightarrow_{c}\infty$ while $r/p\rightarrow0,$ almost surely (a.s.),
\begin{equation}
F_{p,T}\left(  \lambda\right)  \Rightarrow W\left(  \lambda
;c/(1+c),2c/(1+c)\right)  ,\label{convergence}%
\end{equation}
where $\Rightarrow$ denotes the weak convergence of d.f.'s (see Billingsley
(1995), p.191), and $W\left(  \lambda;\gamma_{1},\gamma_{2}\right)  $ denotes
the \textit{Wachter }d.f. with parameters $\gamma_{1}$ and $\gamma_{2}$. The
\textit{Wachter distribution} was derived by Wachter (1980) as the limit of
the empirical distribution of the eigenvalues of the multivariate beta matrix
of growing dimension and degrees of freedom. It has a simple density, which is
introduced in the next section, and, for $\gamma_{2}>\gamma_{1}$ and/or
$\gamma_{2}<1-\gamma_{1},$ point masses at zero and/or one, respectively.

The a.s. weak convergence (\ref{convergence}) and the fact that the squared
sample canonical correlations are no larger than unity imply the a.s.
convergence of averages $\frac{1}{p}%
{\displaystyle\sum_{i=1}^{p}}
f\left(  \lambda_{i}\right)  $ for any $f$ which is bounded and continuous on
$\left[  0,1\right]  $. By definition, the likelihood ratio statistic scaled
by $1/(pT)$ has this form (with omitted first $r$ summands), where
$f(\lambda)=-\log\left(  1-\lambda\right)  $ is continuous but unbounded
function. Therefore, (\ref{convergence}) can guarantee an a.s. asymptotic
lower bound for the scaled LR statistic. For the LR statistic scaled by
$1/p^{2},$ we have, almost surely,%
\begin{equation}
\lim_{p,T\rightarrow_{c}\infty}\inf LR_{r,p,T}/p^{2}\geq-\frac{1}{c}\int%
\log\left(  1-\lambda\right)  \mathrm{d}W\left(  \lambda
;c/(1+c),2c/(1+c)\right)  .\label{bound}%
\end{equation}

In contrast, we show that, under the (standard) asymptotic regime
where\linebreak$T\rightarrow\infty$ while $p$ is held fixed$,$ $LR_{r,p,T}%
/p^{2}$ concentrates around $2$ for relatively large $p$.\footnote{Similar to
(\ref{bound}), our weak convergence results only guarantee that $2 $ is a
lower bound, but we conjecture that it is also the limit of the scaled LR
statistic as first $T\rightarrow\infty$ and then $p\rightarrow\infty$. This
conjecture is supported by Monte Carlo evidence.} A direct calculation reveals
that $2$ is smaller than the lower bound (\ref{bound}), for all $c>0$, with
the gap growing as $c$ increases. That is, the standard asymptotic
distribution of the LR statistic is centered at a too low level, especially
for relatively large $p$. This explains the tendency of the asymptotic LR test
to over-reject the null.

The reason for the poor centering delivered by the standard asymptotic
approximation is that it classifies terms $\left(  p/T\right)  ^{j}$ in the
asymptotic expansion of the likelihood ratio statistic as $O\left(
T^{-j}\right)  .$ When $p$ is relatively large, such terms can substantially
contribute to the finite sample distribution of the statistic, but will be
ignored as asymptotically negligible. In contrast, the \textit{simultaneous}
\textit{asymptotics} classifies all terms $\left(  p/T\right)  ^{j}$ as
$O(1).$ They are not ignored asymptotically, which improves the centering of
the simultaneous asymptotic approximation relative to the standard one.

It is possible to use bound (\ref{bound}), with $c$ replaced by $p/T$, to
construct a Bartlett-type correction factor for the standard LR test. As we
show below, for $p/T<1/3,$ the value of such a theoretical correction factor
is very close to the simulation-based factor described in Johansen, Hansen and
Fachin (2005). However, for larger $p/T$, the values diverge, which may be
caused by the fact that Johansen, Hansen and Fachin's (2005) simulations do
not consider combinations of $p$ and $T$ with $p/T>1/3,$ and the functional
form that they use to fit the simulated correction factors does not work well
uniformly in $p/T$.

The weak convergence result (\ref{convergence}) can be put to a more direct
use by comparing the quantiles of the empirical distribution of the squared
sample canonical correlations with the quantiles of the limiting Wachter
distribution. Under the null, the former quantiles plotted against the latter
ones should form a 45$^{\circ}$ line, asymptotically. Deviations of such a
Wachter quantile-quantile plot from the line indicate violations of the null.
Creating Wachter plots requires practically no additional computations beyond
those needed to compute the LR statistic, and we propose to use this simple
graphical device for a preliminary analysis of cointegration in large VARs.

To the best of our knowledge, our study is the first to derive the limit of
the empirical d.f. of the squared sample canonical correlations between random
walk $X_{t-1}$ and its innovations $\Delta X_{t}$. Wachter (1980) shows that
$W\left(  \lambda;\gamma_{1},\gamma_{2}\right)  $ is the weak limit of the
empirical d.f. of the squared sample canonical correlations between $q$- and
$m$-dimensional independent Gaussian white noises with the size of the sample
$n,$ when $q,m,n\rightarrow\infty$ so that $q/n\rightarrow\gamma_{1}$ and
$m/n\rightarrow\gamma_{2}$. Yang and Pan (2012) show that Wachter's (1980)
result holds without the Gaussianity assumption for i.i.d. data with finite
second moments.

Our proofs do not rely on those previous results. The values of parameters
$\gamma_{1}$ and $\gamma_{2}$ in (\ref{convergence}) imply that the limiting
d.f. for the case of $T$ observations of $p$-dimensional random walk and its
innovations, that we consider in this paper, is the same as the limiting d.f.
for the case of $T+p$ observations of two independent white noises - one
$p$-dimensional and the other $2p$-dimensional. It is tempting to think that
there exists a deep connection between the two cases, even though we were
unable to uncover it so far.

Our paper opens up a new direction for the asymptotic analysis of panel VAR
cointegration tests based on the sample canonical correlations. One such test
is developed in Larsson and Lyhagen (2007). It generalizes Larsson, Lyhagen,
and Lothgren (2001) and Groen and Kleibergen (2003) by allowing for cross-unit
cointegration, which is important from the empirical perspective. Larsson and
Lyhagen (2007) are reluctant to recommend their test for large VARs and
suggest that for the analysis of relatively large panels it may be better to
rely on tighter parameterized models, such as that of Bai and Ng (2004). In
the recent review of the panel cointegration literature, Choi (2015) expresses
a related concern that, with the large number of cross-sectional units,
\textquotedblleft Larsson and Lyhagen's test may not work well even with the
Bartlett's correction.\textquotedblright\ 

We speculate that the Larsson-Lyhagen test, as well as Johansen's LR test,
based on the \textit{simultaneous} asymptotics would work well in panels with
comparable cross-sectional and temporal dimensions. The results of this paper
can be used to describe only the appropriate centering of the corresponding
test statistics. The next step would be to derive the \textit{simultaneous}
asymptotic distribution of scaled deviations of such statistics from the
centering values. We conjecture that the \textit{simultaneous} asymptotic
distribution of $LR_{r,p,T}$ is Gaussian, as is often the case for averages of
regular functions of eigenvalues of large random matrices (see Bai and
Silverstein (2010) and Paul and Aue (2014)). We are currently undertaking work
to validate this conjecture.

The rest of this paper is structured as follows. In Section 2, we prove the
convergence of $F_{p,T}\left(  \lambda\right)  $ to the \textit{Wachter} d.f.
and use this result to derive the asymptotic lower bound for $LR_{r,p,T}.$
Section 3 derives the sequential limit of the empirical d.f. of the squared
sample canonical correlations as, first $T\rightarrow\infty$ and then
$p\rightarrow\infty$. It then uses differences between the obtained sequential
asymptotic limit and the simultaneous limit derived in Section 2 to explain
the over-rejection phenomenon, and to design a theoretical Bartlett-type
correction factor for the LR statistic in high-dimensional VARs. Section 4
contains a Monte Carlo study that confirms good finite sample properties of
the Wachter asymptotic approximation. It also illustrates the proposed Wachter
quantile-quantile plot technique using a relatively high-dimensional
macroeconomic panel. Section 5 concludes and points out directions for future
research. All proofs are given in the Appendix.

\section{Convergence to the Wachter distribution}

Consider the following basic version of (\ref{model})%
\begin{equation}
\Delta X_{t}=\Pi X_{t-1}+\Phi D_{t}+\varepsilon_{t}%
\label{econometricians model}%
\end{equation}
with $d_{D}$-dimensional vector of deterministic regressors $D_{t}$. Let
$R_{0t}$ and $R_{1t}$ be the vectors of residuals from the OLS regressions of
$\Delta X_{t}$ on $D_{t},$ and $X_{t-1}$ on $D_{t},$ respectively. Define%
\begin{equation}
S_{00}=\frac{1}{T}%
{\displaystyle\sum_{t=1}^{T}}
R_{0t}R_{0t}^{\prime},\text{ }S_{01}=\frac{1}{T}%
{\displaystyle\sum_{t=1}^{T}}
R_{0t}R_{1t}^{\prime},\text{ and }S_{11}=\frac{1}{T}%
{\displaystyle\sum_{t=1}^{T}}
R_{1t}R_{1t}^{\prime},\label{SSS}%
\end{equation}
and let $\lambda_{1}\geq...\geq\lambda_{p}$ be the eigenvalues of
$S_{01}S_{11}^{-1}S_{01}^{\prime}S_{00}^{-1}.$

The main goal of this section is to establish the a.s. weak convergence of the
empirical d.f. of the $\lambda$'s to the Wachter d.f., under the null of $r$
cointegrating relationships, when $p,T\rightarrow_{c}\infty$. The Wachter
distribution with d.f. $W\left(  \lambda;\gamma_{1},\gamma_{2}\right)  $ and
parameters $\gamma_{1},\gamma_{2}\in\left(  0,1\right)  $ has density%
\begin{equation}
f_{W}\left(  \lambda;\gamma_{1},\gamma_{2}\right)  =\frac{1}{2\pi\gamma_{1}%
}\frac{\sqrt{\left(  b_{+}-\lambda\right)  \left(  \lambda-b_{-}\right)  }%
}{\lambda\left(  1-\lambda\right)  }\label{densityW}%
\end{equation}
on $\left[  b_{-},b_{+}\right]  \subseteq\left[  0,1\right]  $ with%
\begin{equation}
b_{\pm}=\left(  \sqrt{\gamma_{1}(1-\gamma_{2})}\pm\sqrt{\gamma_{2}%
(1-\gamma_{1})}\right)  ^{2},\label{support boundaries}%
\end{equation}
and atoms of size $\max\left\{  0,1-\gamma_{2}/\gamma_{1}\right\}  $ at zero,
and $\max\left\{  0,1-(1-\gamma_{2})/\gamma_{1}\right\}  $ at unity.

We shall assume that model (\ref{econometricians model}) may be misspecified
in the sense that the true data generating process is described by the
following generalization of (\ref{model})%
\begin{equation}
\Delta X_{t}=\Pi X_{t-1}+%
{\displaystyle\sum_{i=1}^{k-1}}
\Gamma_{i}\Delta X_{t-i}+\Psi F_{t}+\varepsilon_{t},\label{general model}%
\end{equation}
where $\varepsilon_{t},$ $t=1,...,T,$ are still i.i.d. $N(0,\Sigma)$ with
arbitrary $\Sigma>0,$ $\operatorname*{rank}\Pi=r,$ but $k$ is not necessarily
unity, and $F_{t}$ is a $d_{F}$-dimensional vector of deterministic or
stochastic variables that does not necessarily coincide with $D_{t}$. For
example, some of the components of $F_{t}$ may be common factors not observed
and not modelled by the econometrician. Further, we do not put any
restrictions on the roots of the characteristic polynomial associated with
(\ref{general model}). In particular, explosive behavior and seasonal unit
roots are allowed. Finally, no constraints on $F_{t},$ and the initial values
$X_{1-k},...,X_{0}$, apart from the asymptotic requirements on $d_{F}$ and $k$
as spelled out in the following theorem, are imposed.

\begin{theorem}
\label{main}Suppose that the data are generated by (\ref{general model}), and
let $\lambda_{i},$ $i=1,...,p,$ be the eigenvalues of $S_{01}S_{11}^{-1}%
S_{01}^{\prime}S_{00}^{-1},$ where $S_{ij}$ are as defined in (\ref{SSS}).
Further, let $F_{p,T}\left(  \lambda\right)  $ be the empirical d.f. of the
$\lambda$'s, and let $\Gamma=\left[  \Gamma_{1},...,\Gamma_{k-1}\right]  $. If%
\begin{equation}
\frac{1}{p}\left(  d_{D}+d_{F}+r+k+\operatorname*{rank}\Gamma\right)
\rightarrow0\label{wistles}%
\end{equation}
as $p,T\rightarrow_{c}\infty$ while $p$ remains no larger than $T$, then,
almost surely,%
\begin{equation}
F_{p,T}\left(  \lambda\right)  \Rightarrow W\left(  \lambda
;c/(1+c),2c/(1+c)\right)  .\label{5again}%
\end{equation}

\end{theorem}

Condition (\ref{wistles}) requires the number $d_{D}$ of deterministic
regressors in the econometrician's model (\ref{econometricians model}), the
dimensionality $d_{F}$ of $F_{t}$, the number $r$ of the cointegrating
relationships under the null, the order $k$ of the data generating VAR, and
the dimensionality of the union of the column spaces of the matrix
coefficients on \textquotedblleft further lags\textquotedblright\ in
(\ref{general model}) to be either fixed or growing less than proportionally
to the dimensionality $p$ or, equivalently, to the sample size $T$. This
condition rules out situations where some or all lags which are omitted from
the econometrician's model (\ref{econometricians model}) have full rank
coefficients $\Gamma_{i}$. The simplest special situation where (\ref{wistles}%
) is clearly satisfied corresponds to the pure random walk data $\Delta
X_{t}=\varepsilon_{t}$.

The reason why the limit of the empirical d.f. $F_{p,T}(\lambda)$ does not
change when the data generating process (\ref{general model}) changes so that
(\ref{wistles}) remains true is that the corresponding changes in the matrix
$S_{01}S_{11}^{-1}S_{01}^{\prime}S_{00}^{-1}$ have rank that is less than
proportional to $p$ (and to $T$). By the so-called rank inequality (Theorem
A43 in Bai and Silverstein (2010)), the L\'{e}vy distance between the
empirical d.f. of eigenvalues corresponding to versions of $S_{01}S_{11}%
^{-1}S_{01}^{\prime}S_{00}^{-1}$ that differ by a matrix of rank $R$ is no
larger than $R/p,$ which converges to zero as $p,T\rightarrow_{c}\infty.$
Since the L\'{e}vy distance metrizes the weak convergence (see Billingsley
(1995), problem 14.5), the limiting d.f. is not affected. For further details,
see the proof of Theorem \ref{main} in the Appendix.

\begin{remark}
In standard cases where $D_{t}$ is represented by $\left(  1,t\right)  ,$ it
is customary to impose restrictions on $\Phi$ so that there is no quadratic
trend in $X_{t}$ (see Johansen (1995), ch. 6.2). Then, the LR test of the null
of $r$ cointegrating relationships is based on the eigenvalues of
$S_{01}^{\ast}S_{11}^{\ast-1}S_{01}^{\ast\prime}S_{00}^{\ast-1},$ defined
similarly to $S_{01}S_{11}^{-1}S_{01}^{\prime}S_{00}^{-1}$ by replacing
$X_{t-1}$ with $\left(  X_{t-1}^{\prime},t\right)  ^{\prime}$ and regressing
$\Delta X_{t}$ and $\left(  X_{t-1}^{\prime},t\right)  ^{\prime}$ on constant
only to obtain $R_{0t}$ and $R_{1t}.$ The empirical distribution function of
so modified eigenvalues still converges to $W\left(  \lambda
;c/(1+c),2c/(1+c)\right)  $ because the difference between matrices
$S_{01}^{\ast}S_{11}^{\ast-1}S_{01}^{\ast\prime}S_{00}^{\ast-1}$ and
$S_{01}S_{11}^{-1}S_{01}^{\prime}S_{00}^{-1}$ has small rank.
\end{remark}

\begin{figure}[ptb]%
\centering
\includegraphics[
height=3.301in,
width=4.0776in
]%
{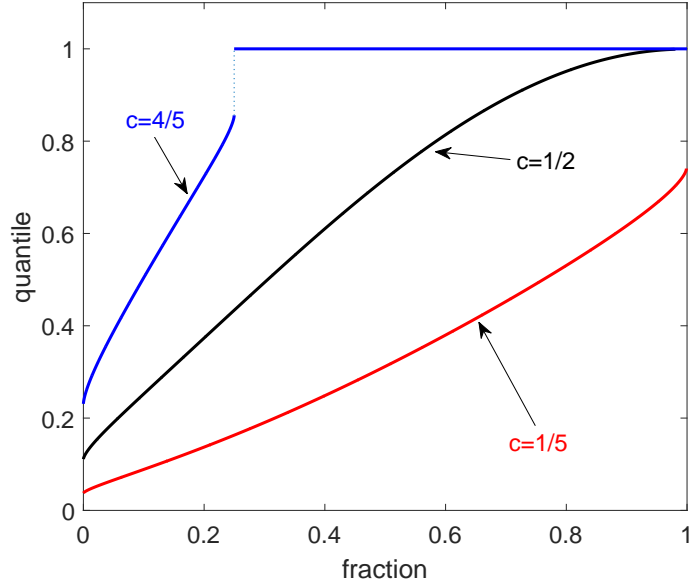}%
\caption{Quantile functions of $W\left(  \lambda;c/\left(  1+c\right)
,2c/\left(  1+c\right)  \right)  $ for $c=1/5,$ $c=1/2,$ and $c=4/5.$}%
\label{illustrationmain}%
\end{figure}

Figure \ref{illustrationmain} shows quantile plots of the Wachter distribution
with parameters $\gamma_{1}=c/(1+c)$ and $\gamma_{2}=2c/(1+c)$ for different
values of $c$. For $c=1/5,$ the dimensionality of the data constitutes 20\% of
the sample size. The corresponding Wachter limit of $F_{p,T}\left(
\lambda\right)  \,$is supported on $\left[  0.04,0.74\right]  $. In
particular, we expect $\lambda_{1}$ be larger than $0.7$ for large $p,T$ even
in the absence of any cointegrating relationships. For $c=1/2,$ the upper
boundary of support of the Wachter limit is unity. This accords with Gonzalo
and Pitarakis' (1995, Lemma 2.3.1) finding that as $T/p\rightarrow2, $
$\lambda_{1}\rightarrow1.$ For $c=4/5,$ the Wachter limit has mass $3/4$ at unity.%

Wachter (1980) derives $W\left(  \lambda;\gamma_{1},\gamma_{2}\right)  $ as
the weak limit of the empirical d.f. of eigenvalues of the $p$-dimensional
beta\footnote{For the definition of the multivariate beta see Muirhead (1982),
p.110.} matrix $B_{p}\left(  n_{1}/2,n_{2}/2\right)  $ with $n_{1},n_{2}$
degrees of freedom as $p,n_{1},n_{2}\rightarrow\infty$ so that $p/n_{1}%
\rightarrow\gamma_{1}/\gamma_{2}$ and $p/n_{2}\rightarrow\gamma_{1}%
/(1-\gamma_{2}). $ The eigenvalues of multivariate beta matrices are related
to many important concepts in multivariate statistics, including canonical
correlations, multiple discriminant ratios, and MANOVA. In particular, the
squared sample canonical correlations between $q$- and $m$-dimensional
independent Gaussian samples of size $n$ are jointly distributed as the
eigenvalues of $B_{q}\left(  m/2,(n-m)/2\right)  ,$ where $q\leq m$ and $n\geq
q+m$. Therefore, their empirical d.f. weakly converges to $W\left(
\lambda;\gamma_{1},\gamma_{2}\right)  $ with $\gamma_{1}=\lim q/n$ and
$\gamma_{2}=\lim m/n,$ as mentioned above. Since the squared canonical
correlations in Theorem \ref{main} are between random walk and its innovations
rather than independent white noises, the convergence to the Wachter
distribution came to us as a pleasant surprise.

In the context of multiple discriminant analysis, Wachter (1976b) proposes to
use a quantile-quantile (qq) plot, where the multiple discriminant ratios are
plotted against quantiles of $W\left(  \lambda;\gamma_{1},\gamma_{2}\right)
$, as a simple graphical method that helps one \textquotedblleft recognize
hopeless from promising analyses at an early stage.\textquotedblright\ A plot
that clearly deviates from the 45$^{\circ}$ line suggests that the data are at
odds with the null hypothesis of the homogeneous population, and a further
analysis of the heterogeneity is useful. Nowadays, such qq plots are called
\textit{Wachter plots} (see Johnstone (2001)).

\begin{figure}[h]%
\centering
\includegraphics[
trim=0.000000in 0.000000in 0.000000in 0.015797in,
height=3.2621in,
width=5.6031in
]%
{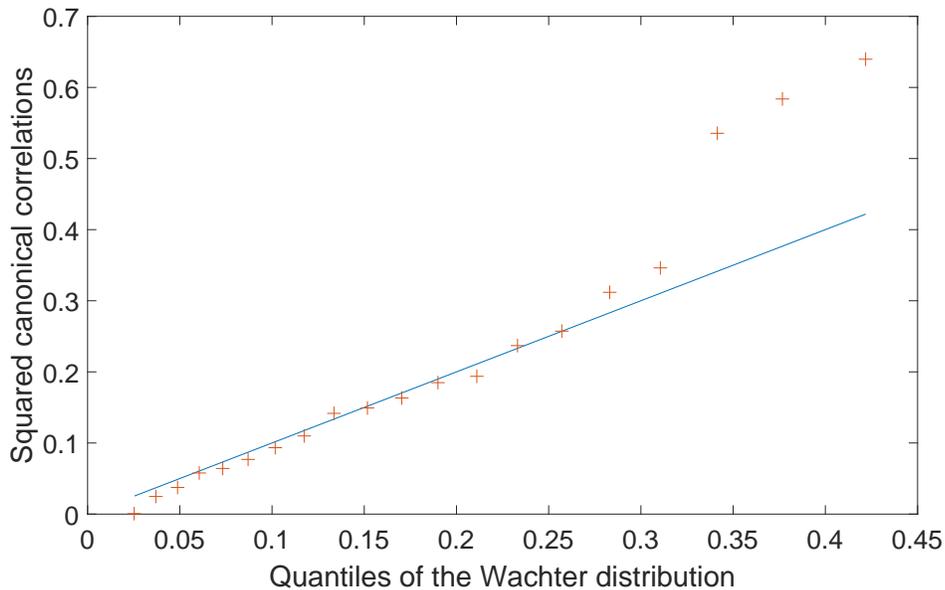}%
\caption{Wachter plot of the squared canonical correlations corresponding to
20-dimensional series with 3 components being white noises and the other
components being independent random walks. $p=20,T=200.$}%
\label{wachterplot}%
\end{figure}

Theorem \ref{main} implies that the Wachter plot can be used as a simple
preliminary assessment of cointegration hypotheses in large VARs. As an
illustration, Figure \ref{wachterplot} shows a Wachter plot of the simulated
sample squared canonical correlations corresponding to a $20$-dimensional
VAR(1) model (\ref{econometricians model}) with $\Pi=\operatorname*{diag}%
\left\{  -I_{3},0\times I_{17}\right\}  $ so that there are three white noise
and seventeen random walk components of $X_{t}$. No deterministic terms are
included. We set $T=200$ so that $c=1/10$. The graph clearly shows three
canonical correlations that destroy the 45$^{\circ}$ line fit, so that the
null hypothesis of no cointegration is compromised.%

\bigskip

Theorem \ref{main} does not provide any explanation to the fact that exactly
three canonical correlations deviate from the 45$^{\circ}$ line in Figure
\ref{wachterplot}. To interpret deviations of the Wachter plots from the
45$^{\circ}$ line, it is desirable to investigate behavior of $F_{p,T}\left(
\lambda\right)  $ under various alternatives. So far, we were able to obtain a
clear result only for the \textquotedblleft extreme\textquotedblright%
\ alternative, where $X_{t}$ is a vector of independent white noises. Under
such an alternative,
\begin{equation}
F_{p,T}\left(  \lambda\right)  \Rightarrow W\left(  \lambda
;c/(2-c),1/(2-c)\right)  .\label{white noise}%
\end{equation}
We plan to publish a full proof of this and some related results elsewhere.

Interestingly, for $c=1/2,$ the Wachter limits (\ref{5again}) and
(\ref{white noise}) corresponding to random walk and white noise nulls,
respectively, coincide. Hence, as $c$ approaches $1/2,$ not only the largest
sample canonical correlation converges to one and the LR test breaks down, but
also the Wachter plot looses the ability to differentiate between opposite
cointegration hypotheses. For smaller values of $c,$ however, the Wachter
limits (\ref{5again}) and (\ref{white noise}) become well separated. We
provide Monte Carlo analysis of the behavior of $F_{p,T}\left(  \lambda
\right)  $ under some alternative hypotheses in Section~4 below.

The almost sure weak convergence of $F_{p,T}\left(  \lambda\right)  $
established in Theorem \ref{main} implies the almost sure convergence of
bounded continuous functionals of $F_{p,T}\left(  \lambda\right)  .$ An
example of such a functional is the scaled Pillai-Bartlett statistic for the
null of no more than $r$ cointegrating relationships (see Gonzalo and
Pitarakis (1995))%
\[
\frac{1}{Tp}PB_{r,p,T}=\frac{1}{p}%
{\displaystyle\sum_{j=r+1}^{p}}
\lambda_{j},
\]
which is asymptotically equivalent to the LR statistic under the standard
asymptotic regime where $p$ is fixed and $T\rightarrow\infty$. Since, by
definition, $\lambda_{j}\in\left[  0,1\right]  ,$ we have%
\begin{equation}
\frac{1}{Tp}PB_{r,p,T}=%
{\displaystyle\int}
f(\lambda)\mathrm{d}F_{p,T}\left(  \lambda\right)  -\frac{1}{p}%
{\displaystyle\sum_{j=1}^{r}}
\lambda_{j},\label{PBrepresentation}%
\end{equation}
where $f$ is the bounded continuous function%
\[
f(\lambda)=\left\{
\begin{array}
[c]{cc}%
0 & \text{for }\lambda<0\\
\lambda & \text{for }\lambda\in\left[  0,1\right] \\
1 & \text{for }\lambda>1.
\end{array}
\right.  .
\]
As long as $r/p\rightarrow0$ as $p,T\rightarrow_{c}\infty,$ the second term on
the right hand side of (\ref{PBrepresentation}) converges to zero. Therefore,
Theorem \ref{main} implies that $PB/\left(  Tp\right)  $ almost surely
converges to $%
{\displaystyle\int}
f(\lambda)\mathrm{d}W\left(  \lambda;c/(1+c),2c/(1+c)\right)  .$ A direct
calculation based on (\ref{densityW}), which we report in the Supplementary
Appendix, yields the following corollary.

\begin{corollary}
\label{PB}Under the assumptions of Theorem \ref{main}, as $p,T\rightarrow
_{c}\infty,$ almost surely,%
\[
PB_{r,p,T}/\left(  Tp\right)  \rightarrow2c/\left(  1+c\right)  +\max\left\{
0,2-1/c\right\}  .
\]

\end{corollary}

A similar analysis of the LR statistic is less straightforward because%
\[
\frac{1}{Tp}LR_{r,p,T}=-\frac{1}{p}%
{\displaystyle\sum_{j=r+1}^{p}}
\log(1-\lambda_{j}),
\]
and $\log\left(  1-\lambda\right)  $ is unbounded on $\lambda\in\left[
0,1\right]  .$ In fact, for $c>1/2,$ $LR_{r,p,T}$ is ill-defined because a
non-negligible proportion of the squared sample canonical correlations exactly
equal unity. However for $c<1/2,$ we can obtain the almost sure asymptotic
lower bound on $LR_{r,p,T}/\left(  Tp\right)  .$ Note that for such $c,$ the
upper bound of the support of $W\left(  \lambda;c/(1+c),2c/(1+c)\right)  $
equals $b_{+}=c\left(  \sqrt{2}-\sqrt{1-c}\right)  ^{-2}<1$. Let
\begin{equation}
\overline{\log}\left(  1-\lambda\right)  =\left\{
\begin{array}
[c]{cc}%
0 & \text{for }\lambda<0\\
\log(1-\lambda) & \text{for }\lambda\in\left[  0,b_{+}\right] \\
\log(1-b_{+}) & \text{for }\lambda>b_{+}.
\end{array}
\right. \label{truncated log}%
\end{equation}
Clearly, $\overline{\log}\left(  1-\lambda\right)  $ is a bounded continuous
function and%
\[
\frac{1}{Tp}LR_{r,p,T}\geq-\frac{1}{p}%
{\displaystyle\sum_{j=r+1}^{p}}
\overline{\log}(1-\lambda_{j}).
\]
Hence, we have the following a.s. lower bound on $LR_{r,p,T}/\left(
Tp\right)  $ (the corresponding calculations are reported in the Supplementary Appendix).

\begin{corollary}
\label{LRsim}Under the assumptions of Theorem \ref{main}, for $c<1/2,$ as
$p,T\rightarrow_{c}\infty,$ almost surely,%
\[
\lim_{p,T\rightarrow_{c}\infty}\inf\frac{1}{Tp}LR_{r,p,T}\geq\frac{1+c}{c}%
\ln\left(  1+c\right)  -\frac{1-c}{c}\ln\left(  1-c\right)  +\frac{1-2c}{c}%
\ln\left(  1-2c\right)  .
\]

\end{corollary}

\begin{remark}
\label{conjectureRemark}We conjecture that the lower bound reported in the
corollary is, in fact, the a.s. limit of $LR_{r,p,T}/\left(  Tp\right)  .$ To
prove this conjecture, one needs to show that $\lambda_{r+1}$ is almost surely
bounded away from unity so that the unboundedness of $\log\left(
1-\lambda\right)  $ is not consequential. We leave this as an important topic
for future research.
\end{remark}

Corollaries \ref{PB} and \ref{LRsim} suggest appropriate \textquotedblleft
centering points\textquotedblright\ for PB and LR statistics scaled by
$1/\left(  Tp\right)  $ for relatively large and comparable $p$ and $T.$ As we
show in the next section, the standard asymptotic distribution of the scaled
PB and LR statistics are likely\footnote{We only establish lower bounds on the
concentration points. However, Monte Carlo evidence suggests that these bounds
are in fact the points of concentration.} to concentrate around very different
points when $p$ becomes large. As will be seen below, this difference sheds
new light on the over-rejection phenomenon discussed above and on the workings
and limitations of the Bartlett correction for the LR statistic. To study the
concentration of the standard asymptotic distributions of the scaled PB and LR
statistics as $p$ grows, we will consider the \textit{sequential} asymptotic
regime, where first $T\rightarrow\infty,$ and then $p\rightarrow\infty$.

\section{Sequential asymptotics and its consequences}

\subsection{Sequential asymptotics}

To obtain useful results under the sequential asymptotics, we shall study
eigenvalues of the scaled matrix%
\begin{equation}
\frac{T}{p}S_{01}S_{11}^{-1}S_{01}^{\prime}S_{00}^{-1}.\label{RawMatrix}%
\end{equation}
Note that under the simultaneous asymptotic regime $p,T\rightarrow_{c}\infty,$
the asymptotic behavior of the scaled and unscaled eigenvalues is the same up
to the factor $c^{-1}$. However, as first $T\rightarrow\infty$ while $p$
remains fixed, the unscaled eigenvalues converge to zero, while scaled ones do
not. We shall denote the empirical d.f. of eigenvalues of the scaled matrix as
$F_{p,T}^{(s)}\left(  \lambda\right)  $.

Without loss of generality, we focus on the case of simple data generating
process
\begin{equation}
\Delta X_{t}=\varepsilon_{t},\text{ }t=1,...,T,\text{ and }X_{0}%
=0,\label{Random Walk DGP}%
\end{equation}
and on the situation, where the econometrician does not include any
deterministic regressors in his or her model, that is $d_{D}=0$. There is no
loss of generality in such simplifications because, as follows from Lemma
\ref{misspec} and the rank inequality used in the proof of Lemma \ref{rank} in
the Appendix, the L\'{e}vy distance between the versions of $F_{p,T}%
^{(s)}\left(  \lambda\right)  $ that correspond to the simplified and the
general cases is bounded from above by a fixed multiple of $\left(
d_{D}+d_{F}+r+k+\operatorname*{rank}\Gamma\right)  /p$. We shall assume that
the latter expression goes to zero as $p\rightarrow\infty.$ Therefore,
whatever the sequential asymptotic limit of $F_{p,T}^{(s)}\left(
\lambda\right)  $ is under the above simplification, it must also be the
sequential asymptotic limit under the general case. For simplicity, in the
rest of this section, we shall assume that $r=0,$ and will consider statistics
$LR_{0,p,T} $ rather than more general $LR_{r,p,T}.$

Under the above simplifications, Johansen's (1988, 1991) results imply that,
as $T\rightarrow\infty$ while $p$ is held fixed, the eigenvalues of the scaled
matrix (\ref{RawMatrix}) jointly converge in distribution to the eigenvalues
of%
\begin{equation}
\frac{1}{p}\int_{0}^{1}\left(  \mathrm{d}B\right)  B^{\prime}\left(  \int%
_{0}^{1}BB^{\prime}\mathrm{d}u\right)  ^{-1}\int_{0}^{1}B\left(
\mathrm{d}B\right)  ^{\prime},\label{JohansenLimit}%
\end{equation}
where $B$ is a $p$-dimensional Brownian motion. We denote the eigenvalues of
(\ref{JohansenLimit}) as $\lambda_{j}^{(\infty)},$ and their empirical d.f. as
$F_{p,\infty}\left(  \lambda\right)  .$

It is not unreasonable to expect that, as $p\rightarrow\infty$, $F_{p,\infty
}\left(  \lambda\right)  $ becomes close to the limit of the empirical
distribution of eigenvalues of (\ref{RawMatrix}) under a simultaneous, rather
than sequential, asymptotic regime $p,T\rightarrow_{\gamma}\infty,$ where
$\gamma$ is close to zero. We shall denote such a limit as $F_{\gamma}\left(
\lambda\right)  .$ This expectation turns out to be correct in the sense that
the following theorem holds.

\begin{theorem}
\label{Levy} Let $F_{0}\left(  \lambda\right)  $ be the weak limit as
$\gamma\rightarrow0$ of $F_{\gamma}\left(  \lambda\right)  .$ Then, as
$p\rightarrow\infty,$ $F_{p,\infty}\left(  \lambda\right)  $ weakly converges
to $F_{0}\left(  \lambda\right)  ,$ in probability.
\end{theorem}

Importantly, the weak limit $F_{0}\left(  \lambda\right)  $ is not the Wachter
d.f. Instead, the following proposition holds.

\begin{proposition}
\label{MPproposition} $F_{0}\left(  \lambda\right)  $ corresponds to a
distribution supported on $\left[  a_{-},a_{+}\right]  $ with
\begin{equation}
a_{\pm}=\left(  1\pm\sqrt{2}\right)  ^{2},\label{MP boundaries}%
\end{equation}
and having density%
\begin{equation}
f\left(  \lambda\right)  =\frac{1}{2\pi}\frac{\sqrt{\left(  a_{+}%
-\lambda\right)  \left(  \lambda-a_{-}\right)  }}{\lambda}.\label{MP density}%
\end{equation}

\end{proposition}

A reader familiar with Large Random Matrix Theory (see Bai and Silverstein
(2010)) might recognize that $F_{0}\left(  \lambda\right)  $ is the cumulative
distribution function of the continuous part of a special case of the
Marchenko-Pastur distribution (Marchenko and Pastur (1967)). The general
Marchenko-Pastur distribution has density%
\[
f_{MP}\left(  \lambda;\kappa,\sigma^{2}\right)  =\frac{1}{2\pi\sigma^{2}%
\kappa}\frac{\sqrt{\left(  a_{+}-\lambda\right)  \left(  \lambda-a_{-}\right)
}}{\lambda}%
\]
over $\left[  a_{-},a_{+}\right]  $ with $a_{\pm}=\sigma^{2}\left(  1\pm
\sqrt{\kappa}\right)  ^{2}$ and a point mass $\max\left\{  0,1-1/\kappa
\right\}  $ at zero. Density (\ref{MP density}) is two times $f_{MP}\left(
\lambda;\kappa,\sigma^{2}\right)  $ with $\kappa=2$ and $\sigma^{2}=1.$ The
multiplication by two is needed because the mass $1/2$ at zero is not a part
of the distribution $F_{0}$.

Recall that, as $T\rightarrow\infty$ while $p$ remains fixed, the LR statistic
converges in distribution to $p$ times the trace of matrix
(\ref{JohansenLimit}):%
\begin{equation}
LR_{0,p,T}\overset{d}{\rightarrow}p%
{\displaystyle\sum_{j=1}^{p}}
\lambda_{j}^{(\infty)}\text{ as }T\rightarrow\infty.\label{Tlimit}%
\end{equation}
On the other hand, according to Theorem \ref{Levy}, for any $\delta_{1}%
,\delta_{2}>0\ $\ and all sufficiently large $p,$%
\begin{equation}
\Pr\left(  \frac{1}{p}%
{\displaystyle\sum_{j=1}^{p}}
\lambda_{j}^{(\infty)}\geq\int\lambda\mathrm{d}F_{0}\left(  \lambda\right)
-\delta_{1}\right)  \geq1-\delta_{2}.\label{plimit}%
\end{equation}
A direct calculation, which we report in the Supplementary Appendix, shows
that $\int\lambda\mathrm{d}F_{0}\left(  \lambda\right)  =2.$ Hence, we have
the following corollary.

\begin{corollary}
\label{Table reduction}As first $T\rightarrow\infty,$ and then $p\rightarrow
\infty,$ the lower probability bound on $LR_{0,p,T}/\left(  2p^{2}\right)  $
is unity in the following sense. As $T\rightarrow\infty$ while $p$ is held
fixed, $LR_{0,p,T}/\left(  2p^{2}\right)  $ converges in distribution to $%
{\displaystyle\sum_{j=1}^{p}}
\lambda_{j}^{(\infty)}/\left(  2p\right)  .$ Further, for any $\delta
_{1},\delta_{2}>0$\ and all sufficiently large $p,$ the probability that $%
{\displaystyle\sum_{j=1}^{p}}
\lambda_{j}^{(\infty)}/\left(  2p\right)  $ is no smaller than $1-\delta_{1}$
is no smaller than $1-\delta_{2}.$
\end{corollary}

The reason why we only claim the lower bound on $LR_{0,p,T}/\left(
2p^{2}\right)  $ is that Theorem \ref{Levy} is silent about the behavior of
the individual eigenvalues $\lambda_{j}^{(\infty)},$ the largest of which may,
in principle, quickly diverge to infinity. We suspect that $2$ is not just the
lower bound, but also the probability limit of $%
{\displaystyle\sum_{j=1}^{p}}
\lambda_{j}^{(\infty)}/p$, so that the sequential probability limit of
$LR_{0,p,T}/\left(  2p^{2}\right)  $ is unity. Verification of this conjecture
requires more work, similar to that discussed in Remark \ref{conjectureRemark}.

Corollary \ref{Table reduction} is consistent with the numerical finding of
Johansen, Hansen and Fachin (2005, Table 2) that, as $T$ becomes large while
$p$ is being fixed, the sample mean of the LR statistic is well approximated
by a polynomial $2p^{2}+\alpha p$ (see also Johansen (1988) and Gonzalo and
Pitarakis (1995)). The value of $\alpha$ depends on how many deterministic
regressors are included in the VAR. Our theoretical result captures only the
`highest order' sequential asymptotic behavior of the LR statistic, which
remains (bounded below by) $2p^{2}$ independent on the number of the
deterministic regressors.

Another piece of numerical support for $2p^{2}$ being not only the lower bound
but also the first order sequential asymptotic approximation to the LR
statistic is provided by the tables of the asymptotic critical values for
Johansen's LR test (see, for example, MacKinnon, Haug and Michelis (1999)).
The critical values in such tables become uncomfortably large for $p>4$. Of
course, the reason for such an unpleasant growth is that those critical values
are of order $2p^{2}$.

The transformation%
\[
LR_{0,p,T}\mapsto LR_{0,p,T}/p-2p
\]
makes the LR statistic `well-behaved' under the sequential asymptotics. The
division by $p$ reduces the `second order behavior' to $O_{\mathrm{P}}(1),$
while subtracting $2p$ eliminates the remaining explosive `highest order
term'. We report the corresponding transformed 95\% critical values alongside
the original ones in Table \ref{CV}.

The transformed critical values resemble 97-99 percentiles of $N(0,1)$. Since
the LR test is one-sided, the resemblance is coincidental. However, we do
expect the sequential asymptotic distribution of the transformed LR statistic
(as well as its simultaneous asymptotic distribution) to be normal (possibly
with non-zero mean and non-unit variance). Our expectation is based on the
fact that $LR_{0,p,T}/p$ behaves as the eigenvalue average (see (\ref{Tlimit}%
)), which is a special case of the so-called linear spectral statistic. The
asymptotic normality of linear spectral statistics for relatively simple
classes of high-dimensional random matrices is a well established result in
the Large Random Matrix Theory (see Bai and Silverstein (2010)). Extending it
to the linear spectral statistics of matrices of form (\ref{JohansenLimit}) is
left as an important direction for future research.%

\begin{table}[tbp] \centering
\begin{tabular}
[c]{l|l|l}\hline
$p$ & Unadjusted CV & CV/$p-2p$\\\hline
\multicolumn{1}{c|}{$1$} & \multicolumn{1}{|c|}{$4.13$} &
\multicolumn{1}{|c}{$2.13$}\\
\multicolumn{1}{c|}{$2$} & \multicolumn{1}{|c|}{$12.32$} &
\multicolumn{1}{|c}{$2.16$}\\
\multicolumn{1}{c|}{$3$} & \multicolumn{1}{|c|}{$24.28$} &
\multicolumn{1}{|c}{$2.09$}\\
\multicolumn{1}{c|}{$4$} & \multicolumn{1}{|c|}{$40.17$} &
\multicolumn{1}{|c}{$2.04$}\\
\multicolumn{1}{c|}{$5$} & \multicolumn{1}{|c|}{$60.06$} &
\multicolumn{1}{|c}{$2.01$}\\
\multicolumn{1}{c|}{$6$} & \multicolumn{1}{|c|}{$83.94$} &
\multicolumn{1}{|c}{$1.99$}\\
\multicolumn{1}{c|}{$7$} & \multicolumn{1}{|c|}{$111.79$} &
\multicolumn{1}{|c}{$1.97$}\\
\multicolumn{1}{c|}{$8$} & \multicolumn{1}{|c|}{$143.64$} &
\multicolumn{1}{|c}{$1.96$}\\
\multicolumn{1}{c|}{$9$} & \multicolumn{1}{|c|}{$179.48$} &
\multicolumn{1}{|c}{$1.94$}\\
\multicolumn{1}{c|}{$10$} & \multicolumn{1}{|c|}{$219.38$} &
\multicolumn{1}{|c}{$1.94$}\\
\multicolumn{1}{c|}{$11$} & \multicolumn{1}{|c|}{$263.25$} &
\multicolumn{1}{|c}{$1.93$}\\
\multicolumn{1}{c|}{$12$} & \multicolumn{1}{|c|}{$311.09$} &
\multicolumn{1}{|c}{$1.92$}\\\hline
\end{tabular}
\caption{The 95\% asymptotic critical values (CV)  for Johansen's LR test.
The unadjsuted values are taken from the first column of Table II in
MacKinnon, Haug and Michelis (1999). }\label{CV}%
\end{table}%

\subsection{Over-rejection phenomenon, and the Bartlett correction}

In this subsection, let us assume that the following conjecture holds.

\begin{conjecture}
\label{conjecture}The simultaneous and sequential asymptotic lower bounds for
the scaled LR statistics derived in Corollaries \ref{LRsim} and
\ref{Table reduction} represent the corresponding simultaneous and sequential
asymptotic limits. Specifically, for $c<1/2$,%
\begin{equation}
\lim_{p,T\rightarrow_{c}\infty}\frac{1}{2p^{2}}LR_{0,p,T}=\frac{1+c}{2c^{2}%
}\ln\left(  1+c\right)  -\frac{1-c}{2c^{2}}\ln\left(  1-c\right)  +\frac
{1-2c}{2c^{2}}\ln\left(  1-2c\right)  ,\label{simlim}%
\end{equation}%
\begin{equation}
\operatorname*{plim}_{p\rightarrow\infty}\lim_{T\rightarrow\infty}\frac
{1}{2p^{2}}LR_{0,p,T}=1.\label{seqlim1}%
\end{equation}

\end{conjecture}

\begin{figure}[ptb]%
\centering
\includegraphics[
height=3.0087in,
width=3.7152in
]%
{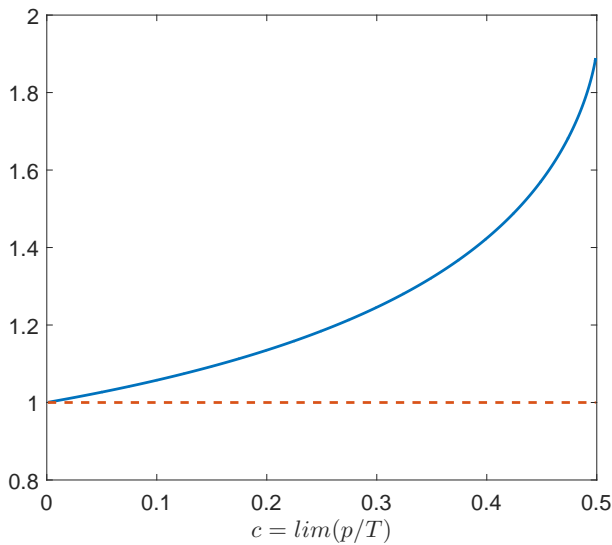}%
\caption{The asymptotic limits (under Conjecture \ref{conjecture}) of the
scaled LR statistic $L_{0,p,T}/\left(  2p^{2}\right)  .$ Dashed line:
sequential asymptotic limit. Solid line: simultaneous asymptotic limit.}%
\label{overrejection}%
\end{figure}

Figure \ref{overrejection} plots the right hand side of (\ref{simlim}) against
the value of $c\in\left[  0,1/2\right)  .$ As demonstrated by the Monte Carlo
analysis of the next section, in finite samples with comparable values of $p$
and $T$, simultaneous asymptotics provides a better approximation to the
finite sample behavior of the LR statistic than the sequential asymptotics.
Therefore, `typical' finite sample values of the LR statistic are concentrated
around the solid line in Figure \ref{overrejection}, and above the dashed
line, which represents the points of concentration of the `standard'
asymptotic critical values for the LR test. In other words, the standard
asymptotic distribution of the LR statistic is centered at a too low level.
This leads to the over-rejection of the null of no cointegration.%

Gonzalo and Pitarakis (1995) propose an interesting approach to address the
problem. Using Monte Carlo, they find that, in contrast to the LR test, the
Pillai-Bartlett test based on the PB statistic\ under-rejects the null.
Therefore, they propose to test cointegration hypotheses using the average of
the LR and PB statistics. According to Corollary \ref{PB}, under the
simultaneous asymptotics $PB/\left(  2p^{2}\right)  \rightarrow1/\left(
1+c\right)  ,$ almost surely. This convergence holds independent on whether
Conjecture \ref{conjecture} is true or not.

\begin{figure}[ptb]%
\centering
\includegraphics[
height=3.0467in,
width=3.7637in
]%
{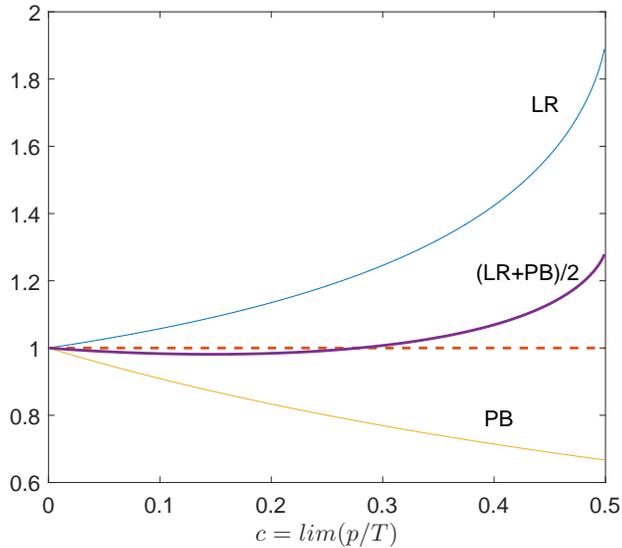}%
\caption{The almost sure limits (under Conjecture \ref{conjecture}) of the
scaled LR, PB, and (LR+PB)/2 statistics under the simulataneous asymptotic
regime.}%
\label{PBaver}%
\end{figure}

The fact that $\left(  1+c\right)  ^{-1}$ is smaller than one, explains the
under-rejection of the test based on the PB statistic. More interestingly, the
average of the simultaneous asymptotic limits of the LR and PB statistics
(divided by $2p^{2}$) turns out to be numerically close to one, and hence to
the point of the concentration of the standard critical values (divided by
$2p^{2}$), at least for $c<1/3.$ Figure \ref{PBaver} shows such an average.
This explains the much better performance of the (LR+PB)/2 test relative to
the LR test in Gonzalo and Pitarakis' (1995) Monte Carlo experiments.%

A more systematic and popular approach to addressing the over-rejection
problem is based on the Bartlett-type correction of the LR statistic. It was
explored in much detail in various important studies, including Johansen
(2002). The idea is to scale the LR statistic so that its finite sample
distribution better fits the asymptotic distribution of the unscaled
statistic. Specifically, let $E_{p,\infty}\left(  LR\right)  $ be the mean of
the asymptotic distribution under the fixed-$p$, large-$T$ asymptotic regime.
Then, if the finite sample mean, $E_{p,T}\left(  LR\right)  $, satisfies%
\begin{equation}
E_{p,T}\left(  LR\right)  =E_{p,\infty}\left(  LR\right)  \left(
1+\frac{a(p)}{T}+o\left(  \frac{1}{T}\right)  \right)
,\label{Bartlett theory}%
\end{equation}
the scaled statistic is defined as $LR/\left(  1+a(p)/T\right)  .$ By
construction, the match between the scaled mean and the original asymptotic
mean is improved by an order of magnitude. Although, as shown by Jensen and
Wood (1997) in the context of unit root testing, the match between higher
moments does not improve by an order of magnitude, it may become substantially
better (see Nielsen (1997)).

A theoretical analysis of the adjustment factor $1+a(p)/T$ can be rather
involved. In general, $a(p)$ will depend not only on $p$, but also on all the
parameters of the VAR. However, for Gaussian VAR(1) without deterministic
terms, under the null of no cointegration, $a(p)$ depends only on $p$.

For $p=1,$ the exact expression for $a(p)$ was derived in Larsson (1998).
Given the difficulty of the theoretical analysis of $a(p),$ Johansen (2002)
proposes to numerically evaluate the Bartlett correction factor $BC_{p,T}%
\equiv E_{p,T}\left(  LR\right)  /E_{p,\infty}\left(  LR\right)  $ by
simulation. Johansen, Hansen and Fachin (2005) simulate $BC_{p,T}$ for various
values of $p\leq10$ and $T\leq3000$ and fit a function of the form%
\[
BC_{p,T}^{\ast}=\exp\left\{  a_{1}\frac{p}{T}+a_{2}\left(  \frac{p}{T}\right)
^{2}+\frac{1}{T}\left[  a_{3}\left(  \frac{p}{T}\right)  ^{2}+b\right]
\right\}
\]
to the obtained results. For relatively large values of $T,$ the term
$\frac{1}{T}\left[  a_{3}\left(  \frac{p}{T}\right)  ^{2}+b\right]  $ in the
above expression is small. When it is ignored, the fitted function becomes
particularly simple:%
\[
\widetilde{BC}_{p,T}=\exp\left\{  0.549\frac{p}{T}+0.552\left(  \frac{p}%
{T}\right)  ^{2}\right\}  .
\]

Our simultaneous and sequential asymptotic results shed light on the workings
of $\widetilde{BC}_{p,T}.$ Given that Conjecture \ref{conjecture} holds,%
\[
\frac{\lim_{p,T\rightarrow_{c}\infty}LR_{0,p,T}}{p\lim_{T\rightarrow
\infty,p\rightarrow\infty}LR_{0,p,T}}=\frac{1+c}{2c^{2}}\ln\left(  1+c\right)
-\frac{1-c}{2c^{2}}\ln\left(  1-c\right)  +\frac{1-2c}{2c^{2}}\ln\left(
1-2c\right)  .
\]
Therefore, for non-negligible $p/T,$ we expect $BC_{p,T}$ to be well
approximated by
\[
\widehat{BC}_{p,T}=\frac{1+\hat{c}}{2\hat{c}^{2}}\ln\left(  1+\hat{c}\right)
-\frac{1-\hat{c}}{2\hat{c}^{2}}\ln\left(  1-\hat{c}\right)  +\frac{1-2\hat{c}%
}{2\hat{c}^{2}}\ln\left(  1-2\hat{c}\right)  ,
\]
where $\hat{c}=p/T$ is the finite sample analog of $c.$

\begin{figure}[ptb]%
\centering
\includegraphics[
height=3.0415in,
width=3.7567in
]%
{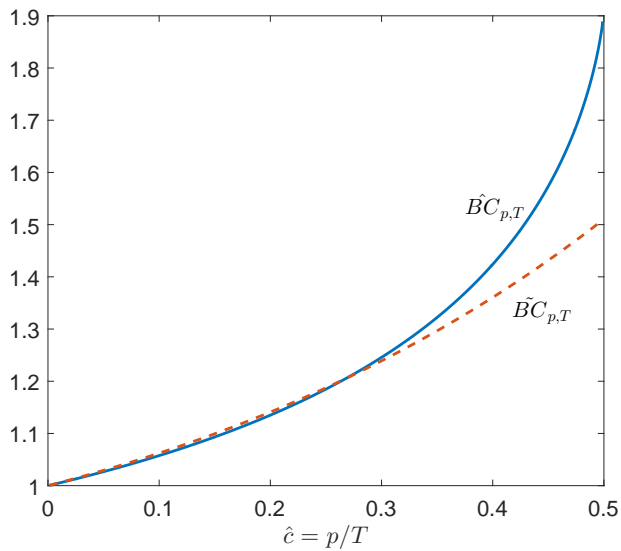}%
\caption{Bartlett correction factors as functions of $p/T.$ Solid line: the
factor based on simultaneous asymptotics. Dashed line: numerical approximation
from Johansen, Hansen and Fachin (2005).}%
\label{bartlettnew}%
\end{figure}

Figure \ref{bartlettnew} superimposes the graphs of $\widehat{BC}_{p,T}$ and
$\widetilde{BC}_{p,T}$ as functions of $\hat{c}.$ For $p/T\leq0.3,$ there is a
strikingly good match between the two curves, with the maximum distance
between them $0.0067$. For $p/T>0.3$ the quality of the match quickly
deteriorates. This can be explained by the fact that all $p,T$-pairs used in
Johansen, Hansen and Fachin's (2005) simulations are such that $p/T<0.3$.%

Further, the good match between $\widetilde{BC}_{p,T}$ and $\widehat{BC}%
_{p,T}$ observed for $p/T<0.3$ would be impossible had Johansen, Hansen and
Fachin's (2005) specified the Bartlett correction factor as a linear function
of $p/T$. Note that the standard theoretical choice for the Bartlett
correction factor, $1+a(p)/T$ from (\ref{Bartlett theory}), can be viewed as a
linear function of $p/T$ with a slope possibly varying with $p$. This is
obvious when $a(p)/T$ is represented as $\frac{p}{T}\beta(p)$ with
$\beta(p)=a(p)/p$. Figure \ref{bartlettnew} shows that such theoretical
correction factors cannot work well uniformly with respect to $p/T$. Uniformly
good correction factors must include terms $\left(  p/T\right)  ^{j}$ with
$j>1$. Under the fixed-$p,$ large-$T$ asymptotics, such terms are of lower
order than $1/T,$ but under the simultaneous asymptotics, they are of order
$O(1)$.

Although the Bartlett-type correction approach may deliver good results for
high-dimensional systems with carefully chosen correction factor, we believe
that tests based on the simultaneous asymptotics of the appropriately scaled
and centered LR statistic would be preferable for relatively large $p$.

\section{Monte Carlo and some examples}

In this section, we describe results of small-scale Monte Carlo experiments
that assess the finite sample quality of the Wachter asymptotic approximation.
In addition, we illustrate the Wachter qq plot technique using a macroeconomic
dataset of relatively high dimensions.

\subsection{Monte Carlo experiments}

\begin{figure}[ptb]%
\centering
\includegraphics[
height=2.5417in,
width=5.61in
]%
{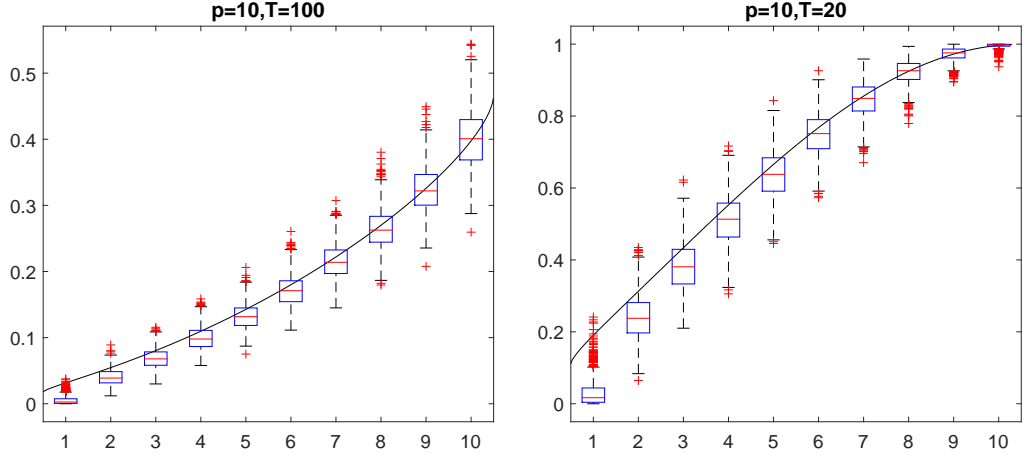}%
\caption{The Tukey boxplots for 1000 MC simulations of ten sample squared
canonical correlations correponding to pure random walk data. The boxplots are
superimposed with the quantile function of the Wachter limit.}%
\label{mcboxplot}%
\end{figure}

First, we generate pure random walk data with zero starting values for
$p=10,T=100$ and $p=10,T=20.$ Throughout this section, the analysis is based
on 1000 Monte Carlo replications. The generated random walk data are
ten-dimensional so that there are ten corresponding squared sample canonical
correlations, $\lambda$. Figure \ref{mcboxplot} shows the Tukey boxplots
summarizing the MC distribution of each of the $\lambda_{i},$ $i=1,...,10$
(sorted in the ascending order throughout this section)$.$ The boxplots are
superimposed with the quantile function of the Wachter limit with $c=1/10$ for
the left panel and $c=1/2$ for the right panel. Precisely, for $x=i,$ we show
the value the $100\left(  i-1/2\right)  /p$ quantile of the Wachter limit. For
$i=1,2,...,10,$ these are the 5-th,15-th,...,95-th quantiles of $W\left(
\lambda;c/\left(  1+c\right)  ,2c/\left(  1+c\right)  \right)  .$

Even for such small values of $p$ and $T,$ the theoretical quantiles track the
location of the MC distribution of the empirical quantiles very well. The
smallest sample canonical correlation is an exception. Its distribution lies
mostly below the corresponding theoretical quantile.

\begin{figure}[ptb]%
\centering
\includegraphics[
height=5.3696in,
width=5.9084in
]%
{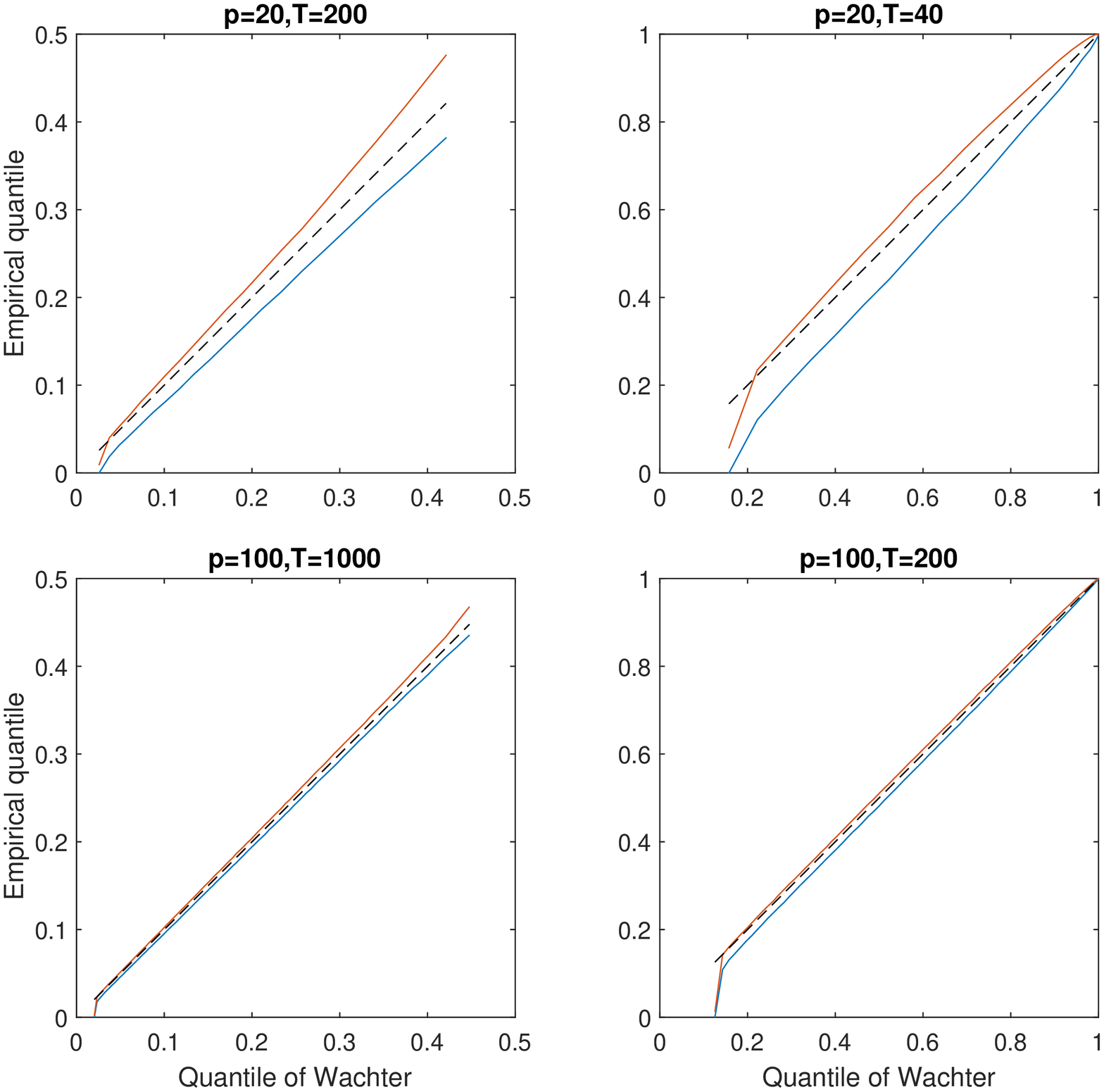}%
\caption{The qq Wachter plots for pure random walk data. The dashed line is
the 45$^{\circ}$ line. The solid lines are the 5-th and the 95-th percentiles
of the MC distributions of $\lambda_{i},$ which are plotted against
$100(i-1/2)/p$ quantiles of the Wachter limit }%
\label{mcqqplot}%
\end{figure}

The dispersion of the MC distributions around the theoretical quantile is
quite large for the chosen small values of $p$ and $T.$ To see how such a
dispersion changes when $p$ and $T$ increase while $p/T$ remains fixed, we
generated pure random walk data with $p=20,T=200$ and $p=100,T=1000$ for
$p/T=1/10$, and with $p=20,T=40$ and $p=100,T=200$ for $p/T=1/2.$ Instead of
reporting the Tukey boxplots, we plot only the 5-th and 95-th percentiles of
the MC distributions of the $\lambda_{i},$ $i=1,...,p$ against $100\left(
i-1/2\right)  /p$ quantiles of the corresponding Wachter limit. The plots are
shown on Figure \ref{mcqqplot}.%

%

We see that the [5\%,95\%] ranges of the MC distributions of $\lambda_{i}$ are
still considerably large for $p=20.$ These ranges become much smaller for
$p=100.$ Interestingly, the distribution of $\lambda_{1}$ remains below the
Wachter limit even for $p=100.$ This does not contradict our theoretical
results because a weak limit of the empirical distribution of $\lambda$'s is
not affected by an arbitrary change in a finite (or slowly growing) number of
them. In fact, we find it somewhat surprising that only the distribution of
$\lambda_{1}$ is not well-alligned with the derived theoretical limit. Our
proofs are based on several low rank alterations of the matrix $S_{01}%
S_{11}^{-1}S_{01}^{\prime}S_{00}^{-1}$, and there is nothing in them that
guarantee that only one eigenvalue of $S_{01}S_{11}^{-1}S_{01}^{\prime}%
S_{00}^{-1}$ behaves in a \textquotedblleft special\textquotedblright\ way. In
future work, it would be interesting to investigate the behavior of
$\lambda_{1}$ and other extreme eigenvalues of $S_{01}S_{11}^{-1}%
S_{01}^{\prime}S_{00}^{-1}$ theoretically.

\begin{figure}[ptb]%
\centering
\includegraphics[
height=5.4034in,
width=5.7674in
]%
{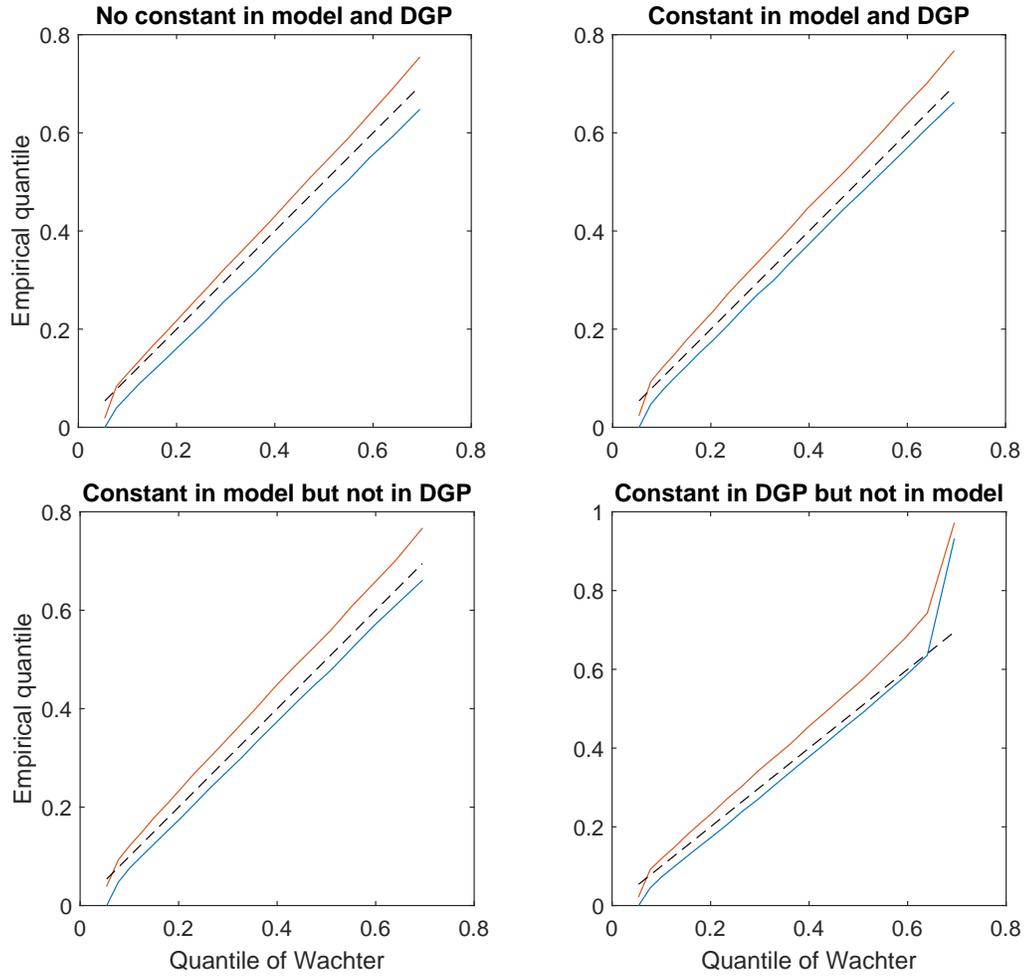}%
\caption{The qq Wachter plots for $p=20$ and $T=100.$ The data generating
process (DGP) is (\ref{general model}) with $k=1,\Pi=0,$ and either $F_{t}=1$
(constant in DGP) or $F_{t}=0 $ (no constant in GDP). The econometrician's
model is (\ref{econometricians model}) with $\Pi=0$ and either $D_{t}=1$
(constant in model) or $D_{t}=0$ (no constant in model).}%
\label{mcmisspec}%
\end{figure}

Next, we explore the effect of the deterministic regressor on the quality of
the Wachter approximation. We generate data with and without constant in the
data generating process (\ref{general model}). That is, we consider two cases:
$F_{t}=1$ and $F_{t}=0.$ The coefficient $\Psi$ on $F_{t}$ is a $N(0,I_{p})$
vector independent across different MC replications. We also consider two
models (\ref{econometricians model}) contemplated by the econometrician: one
with $D_{t}=1$, and the other with $D_{t}=0.$ If $F_{t}\neq D_{t},$ the
econometrican's model is misspecified. Figure \ref{mcmisspec} shows the
Wachter plots similar to those reported in Figure \ref{mcqqplot}. The
dimensions of the data are $p=20$ and $T=100.$

If the data generating process (DGP) contains constant ($F_{t}=1$), but the
econometrician does not include it in his or her model, then the largest
$\lambda,$ $\lambda_{p},$ start to significantly deviate from the 45$^{\circ}$
line on the Wachter plot (lower right panel). If the econometrician's model is
over-specified (lower left panel), there are no dramatic deviations from the line.%

\begin{figure}[ptb]%
\centering
\includegraphics[
height=5.2866in,
width=5.642in
]%
{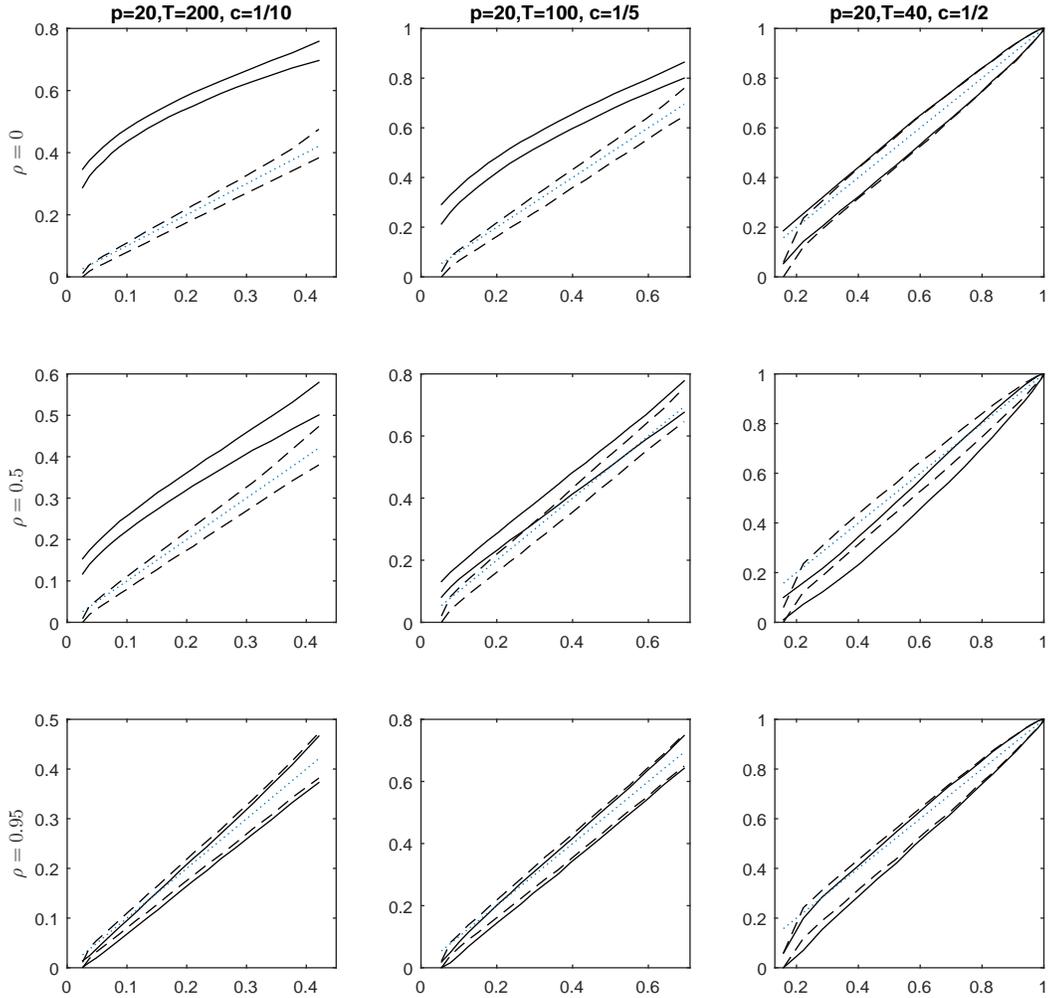}%
\caption{The qq Wachter plots for stationary data $X_{t}=\rho X_{t-1}%
+\varepsilon_{t}.$ Solid lines: 5 and 95 percentiles of the MC distribution of
$\lambda_{i}$ plotted against $100(i-1/2)/p$ quantile of the Wachter limit.
Dashed lines correspond to 5 and 95 percentiles of the MC distribution of
$\lambda_{i}$ for pure random walk data (the null).}%
\label{mcalternatives}%
\end{figure}

Our next Monte Carlo experiment simulates data that are not random walk.
Instead, the data are stationary VAR(1) with zero mean, zero initial value,
and $\Pi=\rho I_{p}.$ We consider three cases of $\rho:0,$ $0.5,$ and $0.95.$
Figure \ref{mcalternatives} shows the Wachter plots with solid lines
representing 5th and 95-th percentiles of the MC distributions of $\lambda
_{i}$ plotted against the $100(i-1/2)/p$ quantiles of the corresponding
Wachter limit. The dashed line correspond to the null case where the data are
pure random walk (shown for comparison).%

The lower panel of the figure corresponds to the most persistent alternative
with $\rho=0.95.$ Samples with $p=20$ seem to be too small to generate
substantial differences in the behavior of Wachter plots under the null and
under such persistent alternatives. The less persistent alternative with
$\rho=0.5$ is easily discriminated against by the Wachter plot for $p/T=1/10$
(left panel). The discrimination power of the plot for $p/T=1/5$ (central
panel) is weaker. For $p/T=1/2$ there is still some discrimination power left,
but the location of the Wachter plot under alternative \textquotedblleft
switches\textquotedblright\ the side relative to the 45$^{\circ}$ line.

The plots easily discriminate against white noise ($\rho=0$) alternative for
$c=1/10$ and $c=1/5,$ but not for $c=1/2.$ In accordance to the result that we
announced above, and plan to publish elsewhere, the Wachter limit for $c=1/2$
approximates equally well the empirical distribution of the squared sample
canonical correlations based on random walk and on white noise data.

Results reported in Figure \ref{mcalternatives} indicate that for relatively
small $p$ and $p/T,$ Wachter plots can be effective in discriminating against
alternatives to the null of no cointegration, where the cointegrating linear
combinations of the data are not very persistent. Further, tests of no
cointegration hypothesis that may be developed using simultaneous asymptotics
would probably need to be two-sided. It is because the location of the Wachter
plot under the alternative may \textquotedblleft switch
sides\textquotedblright\ relative to the 45$^{\circ}$ depending on the
persistence of the data under the alternative. Finally, cases with $c$ close
to 1/2 must be analyzed with much care. For such cases, the behavior of the
sample canonical correlations become similar under extremely different random
walk and white noise data generating processes. Furthermore, the largest
sample canonical correlations are close to unity, which can result in an
unstable behavior of the LR statistic.

\begin{figure}[ptb]%
\centering
\includegraphics[
height=3.096in,
width=5.866in
]%
{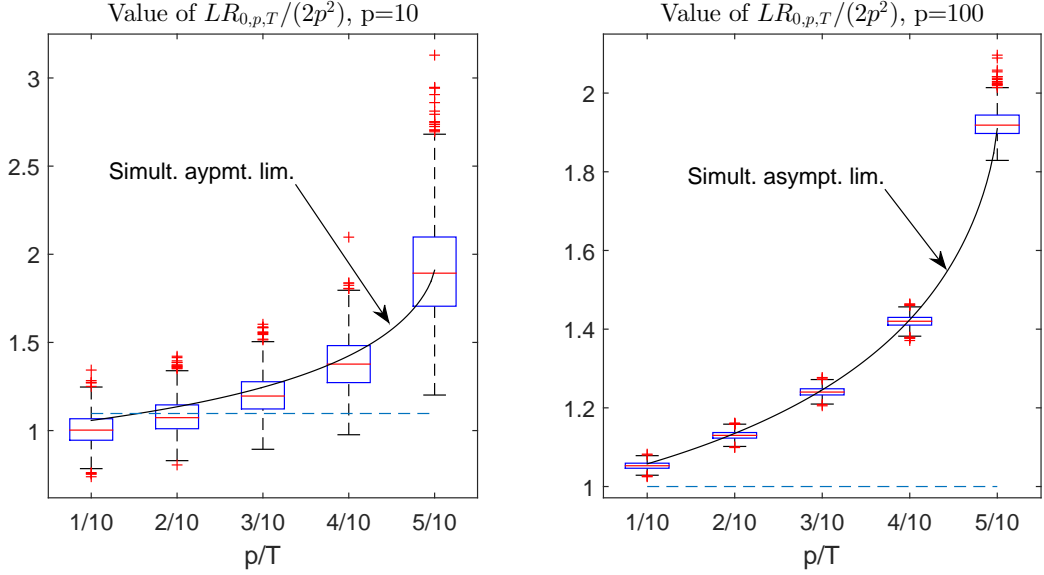}%
\caption{The Tukey boxplots for the MC distributions of $LR_{0,p,T}/\left(
2p^{2}\right)  $ for various $p/T$ ratios. The boxplots are superimposed with
the simultaneous asymptotic lower bound on $LR_{0,p,T}/\left(  2p^{2}\right)
.$ Dashed line in the left panel correspond to 95\% critical value for the
satandard asymptotic Johansen trace test (taken from MacKinnon et al (1999,
Table II)). Dashed line in the right panel has ordinate equal one.}%
\label{mclr}%
\end{figure}

Our final MC experiment studies the finite sample behavior of the scaled LR
statistic $LR_{0,p,T}/\left(  2p^{2}\right)  .$ We simulate pure random walk
data with $p=10$ and $p=100$ and $T$ varying so that $p/T$ equals
1/10,2/10,...,5/10. Corollary \ref{LRsim} shows that the simultaneous
asymptotic lower bound on $LR_{0,p,T}/\left(  2p^{2}\right)  $ has form%
\begin{equation}
\frac{1+c}{2c^{2}}\ln\left(  1+c\right)  -\frac{1-c}{2c^{2}}\ln\left(
1-c\right)  +\frac{1-2c}{2c^{2}}\ln\left(  1-2c\right)
.\label{simultaneous lower bound}%
\end{equation}
Figure \ref{mclr} shows the Tukey boxplots of the MC distributions of
$LR_{0,p,T}/\left(  2p^{2}\right)  $ corresponding to $p/T=1/10,...,5/10$ with
$p=10$ (left panel), and $p=100$ (right panel). The boxplots are superimposed
with the plot of the line representing the above displayed formula for the
lower bound (with $c$ replaced by $p/T$). For the case $p=10,$ we also show
(horizontal dashed line) the 95\% asymptotic critical value (scaled by
$1/(2p^{2})$) of the standard Johansen trace test taken from MacKinnon et al
(1999, Table II). For $p=100,$ critical values for the standard test are not
available, and we show the dashed horizontal line at unit height instead. This
is the sequential asymptotic lower bound on $LR_{0,p,T}/\left(  2p^{2}\right)
$ as established in Corollary \ref{Table reduction}.%

The reported results support our conjecture that the simultaneous asymptotic
lower bound (\ref{simultaneous lower bound}) is, in fact, the simultaneous
asymptotic limit of $LR_{0,p,T}/\left(  2p^{2}\right)  $ for $c<1/2.$
Interestingly, the bound is located near the \textquotedblleft
center\textquotedblright\ of the MC distribution of the scaled LR statistic
even for the case $c=1/2.$

The left panel of Figure \ref{mclr} illustrates the \textquotedblleft
over-rejection phenomenon\textquotedblright. The horizontal dashed line that
corresponds to the 95\% critical value of the standard test is just above the
interquartile range of the MC distribution of $LR_{0,p,T}/\left(
2p^{2}\right)  $ for $c=1/10,$ is below this range for $c\geq3/10,$ and is
below all 1000 MC replications of the scaled LR statistic for $c=5/10$.

Although the lower bound (\ref{simultaneous lower bound}) seems to provide a
very good centering point for the scaled LR statistic, the MC distribution of
this statistic is quite dispersed around such a center for $p=10.$ As
discussed above, we suspect that the scaled statistic centered by
(\ref{simultaneous lower bound}) and appropriately rescaled has Gaussian
simultaneous asymptotic distribution. Optimistically, the Tukey plots on
Figure \ref{mclr}, that correspond to $c<1/2$, look reasonably symmetric
although some skewness is present for the left panel where $p=10.$

\subsection{Examples}

Our first example uses $T=103$ quarterly observations (1973q2-1998q4, with the
initial observation 1973q1) on bilateral US dollar log nominal exchange rates
for $p=17$ OECD countries: Australia, Austria, Belgium, Canada, Denmark,
Finland, France, Germany, Japan, Italy, Korea, Netherlands, Norway, Spain,
Sweden, Switzerland, and the United Kingdom. The data are as in Engel, Mark,
and West (2015), and were downloaded from Charles Engel's website at http: //
www.ssc.wisc.edu / \symbol{126}cengel /. That data are available for a longer
time period up to 2008q1, but we have chosen to use only the \textquotedblleft
early sample\textquotedblright\ that does not include the Euro period.

Engel, Mark, and West (2015) point out that log nominal exchange rates are
well modelled by random walk, but may be cointegrated, which can be utilized
to improve individual exchange rate forecasts relative to the random walk
forecast benchmark. They propose to estimate the common stochastic trends in
the exchange rates by extracting a few factors from the panel. In principle,
the number of factors to extract can be determined using Johansen's test for
cointegrating rank, but Engel, Mark, and West (2015) do not exploit this
possibility, referring to Ho and Sorensen (1996) that indicates poor
performance of the test for large $p.$

\begin{figure}[ptb]%
\centering
\includegraphics[
height=3.0848in,
width=4.0499in
]%
{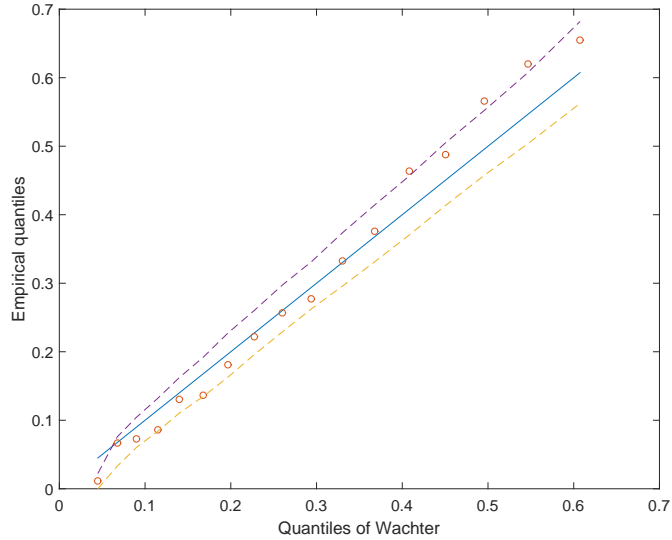}%
\caption{The Wachter plot for the bilateral US log nominal excahnge rates of
17 OECD countries. Dashed lines: 5\% and 95\% quantiles of the MC distribution
of the squared sample canonical correlations under the null of no
cointegration. }%
\label{erfigure}%
\end{figure}

Figure \ref{erfigure} shows the Wachter plot for the log nominal exchange rate
data. The squared sample canonical correlations are computed as the
eigenvalues of $S_{01}S_{11}^{-1}S_{01}^{\prime}S_{00}^{-1},$ where $S_{ij}$
are defined as in (\ref{SSS}) with $R_{0t}$ and $R_{1t}$ being the demeaned
changes and the lagged levels of the log exchange rates, respectively. The
dashed lines correspond to the 5-th and 95-th percentiles of the MC
distribution of the squared canonical correlation coefficients under the null
of no cointegration. Precisely, we generated data from model
(\ref{econometricians model}) with $p=17,$ $T=103$, $\Pi=0,$ $D_{t}=1,$ and
$\Phi$ being i.i.d. $N(0,I_{p})$ vectors across the MC repetitions. Log
exchange rates for 1973q1 was used as the initial value of the generated series.%

The figure shows a mild evidence for cointegration in the data with the
largest five $\lambda$'s being close to the corresponding 95-th percentiles of
the MC distributions. If we interpret this as the existence of five
cointegrating relationships in the data, we would be lead to conclude that
there are twelve stochastic trends. Recall, however, that the ability of the
Wachter plot to differentiate against highly persistent cointegration
alternatives with $p/T\approx1/5$ is very low, so there well may be many more
cointegrating relationships in the data. Whatever such relationships are, the
deviations from the corresponding long-run equilibrium are probably highly
persistent as no dramatic deviations from the 45$^{\circ}$ line are present in
the Wachter plot.

\begin{figure}[ptb]%
\centering
\includegraphics[
height=2.6308in,
width=5.2589in
]%
{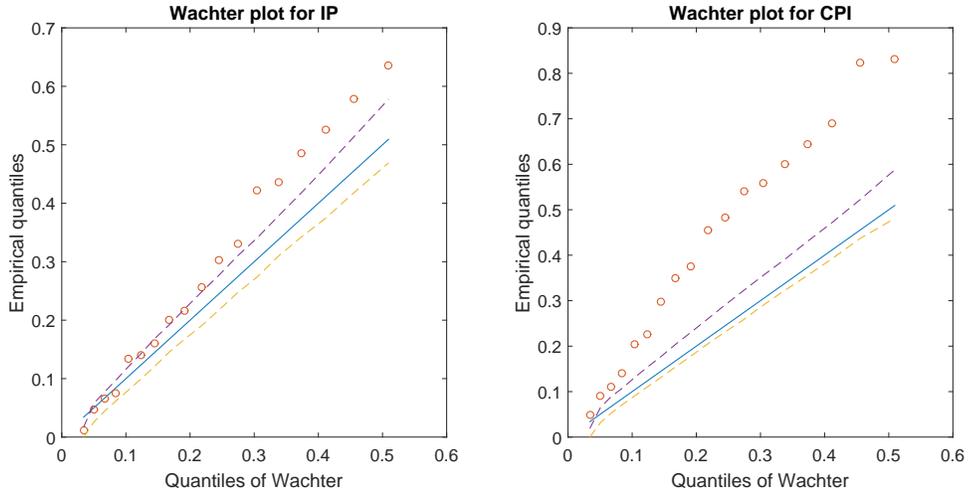}%
\caption{The Wachter plots for the industrial production indices and consumer
price indices of 18 OECD countries. Dashed lines: 5\% and 95\% quantiles of
the MC distribution of the squared sample canonical correlations under the
null of no cointegration. }%
\label{ipcpi}%
\end{figure}

Very different Wachter plots (shown in Figure \ref{ipcpi}) correspond to the
log industrial production (IP) index data and the log consumer price index
(CPI) data for the same countries plus the US. These data are still the same
as in Engel, Mark, and West (2015). We used the long sample 1973q2:2008q1
$(T=140)$ because the IP and CPI data are not affected by the introduction of
the Euro to the same degree as the exchange rate data. For the CPI data, we
included both intercept and trend in model (\ref{econometricians model}) for
the first differences because the level data seem to be quadratically
trending. The plots clearly indicate that the IP and CPI data are either
stationary or cointegrated with potentially many cointegrating relationships,
short run deviations from which are not very persistent.%

\section{Conclusion}

In this paper, we consider the simultaneous, large-$p,$ large-$T$, asymptotic
behavior of the squared sample canonical correlations between $p$-dimensional
random walk and its innovations. We find that the empirical distribution of
these squared sample canonical correlations almost surely weakly converges to
the so-called \textit{Wachter distribution} with parameters that depend only
on the limit of $p/T$ as $p,T\rightarrow_{c}\infty.$ In contrast, under the
sequential asymptotics, when first $T\rightarrow\infty$ and then
$p\rightarrow\infty,$ we establish the convergence in probability to the
so-called Marchenko-Pastur distribution. The differences between the limiting
distributions allow us to explain from a theoretical point of view the
tendency of the LR test for cointegration to severely over-reject the null
when the dimensionality of the data is relatively large. Furthermore, we
derive a simple analytic formula for the Bartlett-type correction factor in
systems with relatively large $p/T$ ratio.

We propose a quick graphical method, the Wachter plot, for a preliminary
analysis of cointegration in large-dimensional systems. The Monte Carlo
analysis shows that the quantiles of the Wachter distribution constitute very
good centering points for the finite sample distributions of the corresponding
squared sample canonical correlations. The quality of the centering is
excellent even for such small $p$ and $T$ as $p=10$ and $T=20.$ However, for
such small values of $p$ and $T,$ the empirical distribution of the squared
sample canonical correlation can considerably fluctuate around the Wachter
limit. As $p$ increases to 100, the fluctuations become numerically very small.

Our analysis leaves many open questions. First, it is very important to study
the fluctuations of the empirical distribution around the Wachter limit. We
conjecture that linear combinations of reasonably smooth functions of the
squared sample canonical correlations, including the $\log(1-\lambda)$ used by
the LR statistic, will be asymptotically Gaussian after appropriate centering
and scaling. The centering can be derived from the results obtained in this
paper. A proof of the asymptotic Gaussianity would require different methods
from those used here. We are currently investigating this research direction.

Further, it would be important to remove the Gaussianity assumption on the
data. We believe that the existence of the finite fourth moments is a
sufficient condition for the validity of the Wachter limit. Next, it would be
interesting to study the simultaneous asymptotic behavior of a few of the
largest sample canonical correlations. This may lead to a modification of
Johansen's maximum eigenvalue test.

Another interesting direction of research is to study situations where the
number of cointegrating relationships under the null is growing proportionally
with $p$ and $T.$ The simultaneous asymptotics of the empirical distribution
of the squared sample canonical correlations under various alternatives, as
well under the null in VAR(k) with $k>1$, also deserves further study.

Still another, totally different, research direction is to investigate the
quality of bootstrap when $p$ is large. Our own very preliminary analysis
indicates that the currently available non-parametric bootstrap procedures
(see, for example, Cavaliere, Rahbek, and Taylor (2012)) do not work well for
$p/T$ as large as, say, 1/3. However, further analysis is needed before we can
claim any specific results. We hope that this paper opens up an interesting
and broad area for future research.

\section{Appendix}

\subsection{Proof of Theorem \ref{main}}

\subsubsection{Reduction to pure random walk data.}

Let $G\left(  \lambda\right)  $ and $\tilde{G}\left(  \lambda\right)  $ be
distribution functions that may depend on $p$ and $T$ and are possibly random.
We shall call them asymptotically equivalent if the a.s. weak convergence
$G\left(  \lambda\right)  \Rightarrow F\left(  \lambda\right)  $ to some
non-random d.f. $F(\lambda)$ implies similar a.s. weak convergence for
$\tilde{G}(\lambda),$ and vice versa. Let $S_{i}$ and $\tilde{S}_{i}$ with
$i=0,1,2$ be, possibly random, matrices that may depend on $p$ and $T$ such
that $S_{i}$ and $\tilde{S}_{i}$ are a.s. positive definite for $i=0,1.$
Below, we shall often refer to the following auxiliary lemma.

\begin{lemma}
\label{rank}If, almost surely, $\frac{1}{p}\operatorname*{rank}\left(
S_{i}-\tilde{S}_{i}\right)  \rightarrow0$ as $p,T\rightarrow_{c}\infty$ for
$i=0,1,2,$ then $G\left(  \lambda\right)  $ and $\tilde{G}(\lambda)$ are
asymptoically equivalent, where $G\left(  \lambda\right)  $ and $\tilde
{G}\left(  \lambda\right)  $ are the empirical d.f. of eigenvalues of
$S_{2}S_{1}^{-1}S_{2}^{\prime}S_{0}^{-1}$ and $\tilde{S}_{2}\tilde{S}_{1}%
^{-1}\tilde{S}_{2}^{\prime}\tilde{S}_{0}^{-1},$ respectively.
\end{lemma}

\paragraph{Proof of Lemma \ref{rank}.}

Let $R=\operatorname*{rank}\left(  S_{2}S_{1}^{-1}S_{2}^{\prime}S_{0}%
^{-1}-\tilde{S}_{2}\tilde{S}_{1}^{-1}\tilde{S}_{2}^{\prime}\tilde{S}_{0}%
^{-1}\right)  .$ The a.s. convergence $\frac{1}{p}\operatorname*{rank}\left(
S_{i}-\tilde{S}_{i}\right)  \rightarrow0$ implies the a.s. convergence
$R/p\rightarrow0.$ On the other hand, by the rank inequality (Theorem A43 in
Bai and Silverstein (2010)), $\mathcal{L}\left(  G,\tilde{G}\right)  \leq
R/p,$ where $\mathcal{L}\left(  G,\tilde{G}\right)  $ is the L\'{e}vy distance
between $G\left(  \lambda\right)  $ and $\tilde{G}(\lambda).$ Recall that the
L\'{e}vy distance metrizes the weak convergence. Therefore, the almost sure
convergence $\mathcal{L}\left(  G,\tilde{G}\right)  \rightarrow0$ yields the
asymptotic equivalence of $G\left(  \lambda\right)  $ and $\tilde{G}%
(\lambda).\square$

Now, let $S_{0}=S_{00},S_{1}=S_{11},$ and $S_{2}=S_{01},$ and let%
\[
\tilde{S}_{0}=\frac{1}{T}%
{\displaystyle\sum_{t=1}^{T}}
\Delta X_{t}\Delta X_{t}^{\prime},\text{ }\tilde{S}_{1}=\frac{1}{T}%
{\displaystyle\sum_{t=1}^{T}}
X_{t-1}X_{t-1}^{\prime},\text{ and }\tilde{S}_{2}=\frac{1}{T}%
{\displaystyle\sum_{t=1}^{T}}
\Delta X_{t}X_{t-1}^{\prime}.
\]
Since $R_{0t}$ and $R_{1t},$ which enter the definition (\ref{SSS}) of
$S_{ij}$, are the residuals in the regressions of $\Delta X_{t}$ on $D_{t}$
and $X_{t-1}$ on $D_{t}$, respectively, we have%
\[
\max_{i=0,1,2}\operatorname*{rank}\left(  S_{i}-\tilde{S}_{i}\right)  \leq
d_{D}.
\]
By assumption, $d_{D}/p\rightarrow0$ as $p,T\rightarrow_{c}\infty,$ so that by
Lemma \ref{rank}, $F_{p,T}(\lambda)$ is asymptotically equivalent to the
empirical d.f. of eigenvalues of $\tilde{S}_{2}\tilde{S}_{1}^{-1}\tilde{S}%
_{2}^{\prime}\tilde{S}_{0}^{-1}.$ Therefore, we may and will replace $R_{0t}$
and $R_{1t}$ in the definitions (\ref{SSS}) of $S_{ij}$ by $\Delta X_{t}$ and
$X_{t-1},$ respectively, without loss of generality. Furthermore, scaling
$S_{ij}$ by $T$ does not change the product $S_{01}S_{11}^{-1}S_{01}^{\prime
}S_{00}^{-1},$ and thus, in the rest of the proof, we shall work with%
\begin{equation}
S_{00}=%
{\displaystyle\sum_{t=1}^{T}}
\Delta X_{t}\Delta X_{t}^{\prime},\text{ }S_{01}=%
{\displaystyle\sum_{t=1}^{T}}
\Delta X_{t}X_{t-1}^{\prime},\text{ and }S_{11}=%
{\displaystyle\sum_{t=1}^{T}}
X_{t-1}X_{t-1}^{\prime}.\label{SSSnew}%
\end{equation}

Next, we show that, still without loss of generality, we may replace the data
generated process (\ref{general model}) by pure random walk with zero initial
value. Indeed, let $X=[X_{-k+1},...,X_{T}],$ where $X_{-k+1},...,X_{0}$ are
arbitrary and $X_{t}$ with $t\geq1$ are generated by (\ref{general model}).
Further, let $\tilde{X}_{-k+1},...,\tilde{X}_{0}$ be zero vectors, $\tilde
{X}_{t}=%
{\displaystyle\sum_{s=1}^{t}}
\varepsilon_{t}$ for $t\geq1,$ and $\tilde{X}=[\tilde{X}_{-k+1},...,\tilde
{X}_{T}].$

\begin{lemma}
\label{misspec}$\operatorname*{rank}\left(  X-\tilde{X}\right)  \leq2\left(
r+\operatorname*{rank}\Gamma+k+d_{F}\right)  .$
\end{lemma}

A proof of this lemma is given in the Supplementary Appendix. It is based on
the representation of $X_{t}$ as a function of the initial values,
$\varepsilon$ and $F$ (see Theorem 2.1 in Johansen (1995)), and requires only
elementary algebraic manipulations. Lemmas \ref{misspec} and \ref{rank}
together with the assumption (\ref{wistles}) imply that replacing $\Delta
X_{t}$ and $X_{t-1}$ in (\ref{SSSnew}) by $\Delta\tilde{X}_{t}$ and $\tilde
{X}_{t-1},$ respectively, does not change the weak limit of $F_{p,T}%
(\lambda).$ Hence, in the rest of the proof of Theorem \ref{main}, without
loss of generality, we shall assume that the data are generated by%
\begin{equation}
\Delta X_{t}=\varepsilon_{t},\text{ }t=1,...,T,\text{ with }X_{0}%
=0.\label{RW DGP}%
\end{equation}

\subsubsection{Block-diagonalization}

Assuming that $\lambda$'s are the eigenvalues of $S_{01}S_{11}^{-1}%
S_{01}^{\prime}S_{00}^{-1}$ with $S_{ij}$ satisfying (\ref{SSSnew}) and
(\ref{RW DGP}), we can interpret them as the squared sample canonical
correlations between lagged values of a random walk $X_{t-1}$ and its current
innovations $\varepsilon_{t}$. Since the sample canonical correlations are
invariant with respect to the multiplication of the data by any invertible
matrix, we assume without loss of generality that the variance of
$\varepsilon_{t}$ equals $\Sigma=I_{p}/T.$ Further, we assume that $T$ is
even. The case of odd $T$ can be analyzed similarly, and we omit it to save space.

Let $\varepsilon=[\varepsilon_{1},...,\varepsilon_{T}]$ and let $U$ be the
upper-triangular matrix with ones above the main diagonal and zeros on the
diagonal. Then $\varepsilon U=\left[  X_{0},...,X_{T-1}\right]  ,$ so that
\begin{equation}
S_{00}=\varepsilon\varepsilon^{\prime},\text{ }S_{01}=\varepsilon U^{\prime
}\varepsilon^{\prime},\text{ and }S_{11}=\varepsilon UU^{\prime}%
\varepsilon^{\prime}.\label{SSSu}%
\end{equation}
We shall show that the empirical d.f. of the $\lambda$'s, $F_{p,T}\left(
\lambda\right)  ,$ is asymptotically equivalent to the empirical d.f. $\hat
{F}_{p,T}\left(  \lambda\right)  $ of eigenvalues of $CD^{-1}C^{\prime}%
A^{-1},$ where
\[
C=\varepsilon\Delta_{2}^{\prime}\varepsilon^{\prime},\text{ }D=\varepsilon
\Delta_{1}\varepsilon^{\prime},\text{ and }A=\varepsilon\varepsilon^{\prime},
\]
$\Delta_{1}$ is a diagonal matrix,%
\begin{equation}
\Delta_{1}=\operatorname*{diag}\left\{  r_{1}^{-1}I_{2},...,r_{T/2}^{-1}%
I_{2}\right\}  ,\label{D1}%
\end{equation}
and $\Delta_{2}$ is a block-diagonal matrix,%
\begin{equation}
\Delta_{2}=\operatorname*{diag}\left\{  r_{1}^{-1}\left(  R_{1}-I_{2}\right)
,...,r_{T/2}^{-1}\left(  R_{T/2}-I_{2}\right)  \right\}  .\label{D2}%
\end{equation}
Here $I_{2}$ is the 2-dimensional identity matrix, and $r_{j},R_{j}$ are
defined as follows. Let $\theta=-2\pi/T.$ Then for $j=1,2,...$%
\[
r_{j+1}=2-2\cos j\theta,\text{ }R_{j+1}=\left(
\begin{array}
[c]{cc}%
\cos j\theta & -\sin j\theta\\
\sin j\theta & \cos j\theta
\end{array}
\right)  ,
\]
whereas $r_{1}=4,$ $R_{1}=-I_{2}.$

\begin{lemma}
\label{Lemma 1}The distribution functions $F_{p,T}\left(  \lambda\right)  $
and $\hat{F}_{p,T}\left(  \lambda\right)  $ are asymptotically equivalent.
\end{lemma}

\paragraph{Proof of Lemma \ref{Lemma 1}.}

Let $V$ be the circulant matrix (see Golub and Van Loan (1996, p.201)) with
the first column $v=\left(  -1,1,0,...,0\right)  ^{\prime}.$ Direct
calculations show that $UV=I_{T}-le_{T}^{\prime}$ and $VU=I_{T}-e_{1}%
l^{\prime},$ where $e_{j}$ is the $j$-th column of $I_{T},$ and $l$ is the
vector of ones. Using these identities, it is straightforward to verify that%
\begin{align}
U  & =\left(  V+e_{1}e_{1}^{\prime}\right)  ^{-1}-le_{1}^{\prime},\text{
and}\label{identity1}\\
UU^{\prime}  & =\left(  V^{\prime}V-\left(  e_{1}-e_{T}\right)  \left(
e_{1}-e_{T}\right)  ^{\prime}+e_{T}e_{T}^{\prime}\right)  ^{-1}-ll^{\prime
}.\label{identity2}%
\end{align}
Now, let us define
\[
C_{1}=\varepsilon\left(  U+le_{1}^{\prime}\right)  ^{\prime}\varepsilon
^{\prime}\text{ and }D_{1}=\varepsilon\left(  UU^{\prime}+ll^{\prime}\right)
\varepsilon^{\prime}.
\]
Using identities (\ref{SSSu}) for $S_{ij}$ and Lemma \ref{rank}, we conclude
that $F_{p,T}(\lambda)$ is asymptotically equivalent to $F_{p,T}^{(1)}\left(
\lambda\right)  $, where $F_{p,T}^{(1)}\left(  \lambda\right)  $ is the
empirical d.f. of the eigenvalues of $C_{1}D_{1}^{-1}C_{1}^{\prime}A^{-1}.$

Further, (\ref{identity1}) and (\ref{identity2}) yield%
\begin{align*}
C_{1}  & =\varepsilon\left(  V+e_{1}e_{1}^{\prime}\right)  ^{-1}%
\varepsilon^{\prime}\text{ and }\\
D_{1}  & =\varepsilon\left(  V^{\prime}V-\left(  e_{1}-e_{T}\right)  \left(
e_{1}-e_{T}\right)  ^{\prime}+e_{T}e_{T}^{\prime}\right)  ^{-1}\varepsilon
^{\prime}.
\end{align*}
Applying Lemma \ref{rank} one more time, we obtain the asymptotic equivalence
of $F_{p,T}^{(1)}(\lambda)$ and $F_{p,T}^{(2)}\left(  \lambda\right)  ,$ where
$F_{p,T}^{(2)}\left(  \lambda\right)  $ is the empirical d.f. of the
eigenvalues of $C_{2}D_{2}^{-1}C_{2}^{\prime}A^{-1}$ with%
\begin{equation}
C_{2}=\varepsilon V^{-1}\varepsilon^{\prime}\text{ and }D_{2}=\varepsilon
\left(  V^{\prime}V\right)  ^{-1}\varepsilon^{\prime}.\label{C2D2}%
\end{equation}
As is well known (see, for example, Golub and Van Loan (1996), chapter 4.7.7),
$T\times T$ circulant matrices can be expressed in terms of the discrete
Fourier transform matrices%
\[
\mathcal{F}=\left\{  \exp\left(  \mathrm{i}\theta\left(  s-1\right)  \left(
t-1\right)  \right)  \right\}  _{s,t=1}^{T}%
\]
with $\theta=-2\pi/T.$ Precisely,%
\[
V=\frac{1}{T}\mathcal{F}^{\ast}\operatorname*{diag}\left(  \mathcal{F}%
v\right)  \mathcal{F},\text{ and }V^{\prime}V=\frac{1}{T}\mathcal{F}^{\ast
}\operatorname*{diag}\left(  \mathcal{F}w\right)  \mathcal{F},
\]
where $w=\left(  2,-1,0,...,0,-1\right)  ^{\prime}$ and the star superscript
denotes transposition and complex conjugation. For the $s$-th diagonal
elements of $\operatorname*{diag}\left(  \mathcal{F}v\right)  $ and
$\operatorname*{diag}\left(  \mathcal{F}w\right)  ,$ we have%
\[
\operatorname*{diag}\left(  \mathcal{F}v\right)  _{s}=-1+\exp\left\{
\mathrm{i}\theta\left(  s-1\right)  \right\}  ,\text{ and }%
\operatorname*{diag}\left(  \mathcal{F}w\right)  _{s}=2-2\cos\left(
s-1\right)  \theta.
\]
Note that $\operatorname*{diag}\left(  \mathcal{F}w\right)  _{s}%
=\operatorname*{diag}\left(  \mathcal{F}w\right)  _{T+2-s}$ for $s=2,3,...$ If
$T$ is even, as we assumed above, then there are $T/2-1$ pairs $\left(
s,T+2-s\right)  $, and there is one pair $\left(  1,T/2+1\right)  $ that
correspond to%
\[
\operatorname*{diag}\left(  \mathcal{F}w\right)  _{1}=0,\operatorname*{diag}%
\left(  \mathcal{F}w\right)  _{T/2+1}=4.
\]

Define a permutation matrix $P$ so that the equal diagonal elements of
$P^{\prime}\operatorname*{diag}\left(  \mathcal{F}w\right)  P$ are grouped in
adjacent pairs. Precisely, let $P=\left\{  p_{st}\right\}  $, where%
\[
p_{st}=\left\{
\begin{array}
[c]{l}%
1\text{ if }t=2s-1\text{ for }s=1,...,T/2\\
1\text{ if }t=2\left(  T-s+2\right)  \operatorname{mod}T\text{ for
}s=T/2+1,...,T\\
0\text{ otherwise}%
\end{array}
\right.
\]
and let $W$ be the unitary matrix%
\[
W=\left(
\begin{array}
[c]{cc}%
I_{2} & 0\\
0 & I_{T/2}\otimes Z
\end{array}
\right)  \text{ with }Z=\frac{1}{\sqrt{2}}\left(
\begin{array}
[c]{cc}%
1 & 1\\
\mathrm{i} & -\mathrm{i}%
\end{array}
\right)  ,
\]
where $\otimes$ denotes the Kronecker product. Further, let $Q=\frac{1}%
{\sqrt{T}}WP^{\prime}\mathcal{F}$. As is easy to check, $Q$ is an orthogonal
matrix. Furthermore,%
\[
V=Q^{\prime}\left(  \Delta_{2}^{-1}+2e_{1}e_{1}^{\prime}\right)  Q,\text{ and
}VV^{\prime}=Q^{\prime}\left(  \Delta_{1}^{-1}-4e_{1}e_{1}^{\prime}\right)  Q,
\]
where $\Delta_{1}$ and $\Delta_{2}$ are as defined in (\ref{D1}) and
(\ref{D2}). Combining this with (\ref{C2D2}) and using Lemma \ref{rank} once
again, we obtain the asymptotic equivalence of $F_{p,T}^{(2)}(\lambda)$ and
$F_{p,T}^{(3)}\left(  \lambda\right)  ,$ where $F_{p,T}^{(3)}\left(
\lambda\right)  $ is the empirical d.f. of the eigenvalues of $C_{3}D_{3}%
^{-1}C_{3}^{\prime}A^{-1}$ with%
\[
C_{3}=\varepsilon Q^{\prime}\Delta_{2}Q\varepsilon^{\prime}\text{ and }%
D_{3}=\varepsilon Q^{\prime}\Delta_{1}Q\varepsilon^{\prime}.
\]
Because of the rotational invariance of the Gaussian distribution, the
distributions of $\varepsilon Q^{\prime}$ and $\varepsilon$ are the same.
Hence, $F_{p,T}^{(3)}\left(  \lambda\right)  $ is asymptotically equivalent to
$\hat{F}_{p,T}\left(  \lambda\right)  ,$ and thus, $\hat{F}_{p,T}\left(
\lambda\right)  $ is asymptotically equivalent to $F_{p,T}\left(
\lambda\right)  $.$\square$

\subsubsection{A system of equations for the Stieltjes transform}

Our proof of the almost sure weak convergence of $\hat{F}_{p,T}\left(
\lambda\right)  $ to the Wachter distribution consists of showing that the
Stieltjes transform of $\hat{F}_{p,T}\left(  \lambda\right)  $,%
\begin{equation}
\hat{m}_{p,T}(z)=\int\frac{1}{\lambda-z}\hat{F}_{p,T}\left(  \mathrm{d}%
\lambda\right)  ,\label{mhat def}%
\end{equation}
almost surely converges pointwise in $z\in\mathbb{C}^{+}=\left\{
\zeta:\mathfrak{I}\zeta>0\right\}  $ to the Stieltjes transform $m(z)$ of the
Wachter distribution. To establish such a convergence, we show that, if $m$ is
a limit of $\hat{m}_{p,T}(z)$ along any \textit{subsequence} of
$p,T\rightarrow_{c}\infty,$ then it must satisfy a system of equations with
unique solution given by $m(z)$. The almost sure convergence of $\hat{F}%
_{p,T}\left(  \lambda\right)  $ (and thus, also of $F_{p,T}\left(
\lambda\right)  $) to the Wachter distribution follows then from the
Continuity Theorem for the Stieltjes transforms (see, for example, Corollary 1
in Geronimo and Hill (2003)).

We shall write $\hat{m}$ for the Stieltjes transform $\hat{m}_{p,T}(z)$ to
simplify notation. Let
\begin{equation}
M=CD^{-1}C^{\prime}-zA\text{ and }\tilde{M}=C^{\prime}A^{-1}%
C-zD.\label{MMtilde}%
\end{equation}
Then by definition (\ref{mhat def}), $\hat{m}$ must satisfy the following
equations%
\begin{align}
\hat{m}  & =\frac{1}{p}\operatorname*{tr}\left[  AM^{-1}\right]
,\label{equation2}\\
\hat{m}  & =\frac{1}{p}\operatorname*{tr}\left[  D\tilde{M}^{-1}\right]
.\label{equation4}%
\end{align}

Let us study the above traces in detail. Define%
\[
\varepsilon_{(j)}=\left[  \varepsilon_{2j-1},\varepsilon_{2j}\right]  ,\text{
}j=1,...,T/2.
\]
We now show that the traces in (\ref{equation2}) and (\ref{equation4}) can be
expressed as functions of the terms having form $\varepsilon_{(j)}^{\prime
}\Omega_{j}\varepsilon_{(j)},$ where $\Omega_{j}$ is independent from
$\varepsilon_{(j)}.$ Then, we argue that
\[
\varepsilon_{(j)}^{\prime}\Omega_{j}\varepsilon_{(j)}-\frac{1}{T}%
\operatorname*{tr}\left[  \Omega_{j}\right]  I_{2}%
\]
a.s. converge to zero, and use this fact to derive equations that the limit of
$\hat{m},$ if it exists, must satisfy.

First, consider (\ref{equation2}). Note that%
\begin{equation}
\frac{1}{p}\operatorname*{tr}\left[  AM^{-1}\right]  =\frac{1}{p}%
{\displaystyle\sum_{j=1}^{T/2}}
\operatorname*{tr}\left[  \varepsilon_{(j)}^{\prime}M^{-1}\varepsilon
_{(j)}\right]  .\label{germ of eq2}%
\end{equation}
Let us introduce new notation:%
\[
\Delta_{1j}=r_{j}^{-1}I_{2},\text{ }\Delta_{2j}=r_{j}^{-1}\left(  R_{j}%
-I_{2}\right)  ,
\]%
\[
C_{j}=C-\varepsilon_{(j)}\Delta_{2j}^{\prime}\varepsilon_{(j)}^{\prime},\text{
}D_{j}=D-\varepsilon_{(j)}\Delta_{1j}\varepsilon_{(j)}^{\prime},
\]%
\[
A_{j}=A-\varepsilon_{(j)}\varepsilon_{(j)}^{\prime},\text{ and }M_{j}%
=C_{j}D_{j}^{-1}C_{j}^{\prime}-zA_{j}.
\]
In addition, let%
\begin{align*}
s_{j}  & =\varepsilon_{(j)}^{\prime}D_{j}^{-1}\varepsilon_{(j)},u_{j}%
=\varepsilon_{(j)}^{\prime}D_{j}^{-1}C_{j}^{\prime}M_{j}^{-1}\varepsilon
_{(j)},\\
v_{j}  & =\varepsilon_{(j)}^{\prime}M_{j}^{-1}\varepsilon_{(j)},\text{ and }\\
w_{j}  & =\varepsilon_{(j)}^{\prime}D_{j}^{-1}C_{j}^{\prime}M_{j}^{-1}%
C_{j}D_{j}^{-1}\varepsilon_{(j)}.
\end{align*}

A straightforward algebra that involves multiple use of the
Sherman-Morrison-Woodbury formula (see Golub and Van Loan (1996), p.50)
\begin{equation}
\left(  V+XWY\right)  ^{-1}=V^{-1}-V^{-1}X\left(  W^{-1}+YV^{-1}X\right)
^{-1}YV^{-1},\label{identity}%
\end{equation}
and the identity%
\begin{equation}
\Delta_{2j}\Delta_{2j}^{\prime}=\Delta_{2j}^{\prime}\Delta_{2j}=\Delta
_{1j},\label{delta identities}%
\end{equation}
establishes the following equality%
\begin{equation}
\varepsilon_{(j)}^{\prime}M^{-1}\varepsilon_{(j)}=v_{j}-[v_{j},u_{j}^{\prime
}]\Omega_{j}[v_{j},u_{j}^{\prime}]^{\prime},\label{stem of eq2}%
\end{equation}
where%
\[
\Omega_{j}=\left(
\begin{array}
[c]{cc}%
\frac{1}{1-z}I_{2}+v_{j} & \frac{1}{1-z}r_{j}\Delta_{2j}^{\prime}%
+u_{j}^{\prime}\\
\frac{1}{1-z}r_{j}\Delta_{2j}+u_{j} & \frac{z}{1-z}r_{j}I_{2}-s_{j}+w_{j}%
\end{array}
\right)  ^{-1}.
\]
A derivation of (\ref{stem of eq2}) can be found in the Supplementary Appendix.

Let us define
\begin{align*}
\hat{s}  & =\frac{1}{T}\operatorname*{tr}\left[  D^{-1}\right]  ,\text{ }%
\hat{u}=\frac{1}{T}\operatorname*{tr}\left[  D^{-1}C^{\prime}M^{-1}\right]
,\\
\hat{v}  & =\frac{1}{T}\operatorname*{tr}\left[  M^{-1}\right]  ,\text{ and
}\\
\hat{w}  & =\frac{1}{T}\operatorname*{tr}\left[  D^{-1}C^{\prime}M^{-1}%
CD^{-1}\right]  .
\end{align*}
We have the following lemma, where $\left\Vert \cdot\right\Vert $ denotes the
spectral norm. Its proof is given in the Supplementary Appendix.

\begin{lemma}
\label{Rigour1}For all $z\in\mathbb{C}^{+},$ as $p,T\rightarrow_{c}\infty, $
we have%
\begin{align*}
& \max_{j=1,...,T/2}\left\Vert s_{j}-\hat{s}I_{2}\right\Vert
\overset{a.s}{\rightarrow}0,\text{ }\max_{j=1,...,T/2}\left\Vert u_{j}-\hat
{u}I_{2}\right\Vert \overset{a.s}{\rightarrow}0\\
& \max_{j=1,...,T/2}\left\Vert v_{j}-\hat{v}I_{2}\right\Vert
\overset{a.s}{\rightarrow}0,\text{ }\max_{j=1,...,T/2}\left\Vert w_{j}-\hat
{w}I_{2}\right\Vert \overset{a.s}{\rightarrow}0.
\end{align*}

\end{lemma}

The lemma yields an approximation to the right hand side of (\ref{stem of eq2}%
), which we use in (\ref{germ of eq2}) and (\ref{equation2}) to obtain the
following result.

\begin{proposition}
\label{the second equation}There exists $\zeta>0$ such that, for any $z$ with
zero real part, $\mathfrak{R}z=0$, and the imaginary part satisfying
$\mathfrak{I}z>\zeta,$ we have%
\begin{equation}
\hat{m}=\frac{1}{2\pi c}\int_{0}^{2\pi}\frac{f_{1}\left(  \varphi\right)
}{\left(  1-z\right)  f_{1}\left(  \varphi\right)  +f_{2}\left(
\varphi\right)  }\mathrm{d}\varphi+o(1),\text{ where}\label{asyeq1}%
\end{equation}%
\begin{align*}
f_{1}\left(  \varphi\right)   & =\left(  \hat{w}-\hat{s}-4\sin^{2}%
\varphi\right)  \hat{v}-\hat{u}^{2},\\
f_{2}\left(  \varphi\right)   & =\hat{w}-\hat{s}-4\sin^{2}\varphi\left(
1-\hat{u}-\hat{v}\right)  ,
\end{align*}
and $o(1)\overset{a.s}{\rightarrow}0$, as $p,T\rightarrow_{c}\infty.$
\end{proposition}

\paragraph{Proof of Proposition \ref{the second equation}.}

Consider a $2\times2$ matrix $\hat{S}_{j}$ that is obtained from
$\varepsilon_{(j)}^{\prime}M^{-1}\varepsilon_{(j)}$ by replacing $s_{j}%
,v_{j},u_{j}$ and $w_{j}$ in (\ref{stem of eq2}) with $\hat{s}I_{2},\hat
{v}I_{2},\hat{u}I_{2},$ and $\hat{w}I_{2},$ respectively. We have%
\[
\hat{S}_{j}=\hat{v}I_{2}-[\hat{v}I_{2},\hat{u}I_{2}]\hat{\Omega}_{j}[\hat
{v}I_{2},\hat{u}I_{2}]^{\prime},
\]
where%
\[
\hat{\Omega}_{j}=\left(
\begin{array}
[c]{cc}%
\frac{1}{1-z}I_{2}+\hat{v}I_{2} & \frac{1}{1-z}r_{j}\Delta_{2j}^{\prime}%
+\hat{u}I_{2}\\
\frac{1}{1-z}r_{j}\Delta_{2j}+\hat{u}I_{2} & \frac{z}{1-z}r_{j}I_{2}+(\hat
{w}-\hat{s})I_{2}%
\end{array}
\right)  ^{-1}.
\]
A simple algebra and the identity $\Delta_{2j}+\Delta_{2j}^{\prime}=-I_{2} $
yield%
\begin{align}
\hat{\Omega}_{j}  & =\frac{1-z}{\delta_{j}}\tilde{\Omega}_{j},\text{
where}\label{Omega}\\
\tilde{\Omega}_{j}  & =\left(
\begin{array}
[c]{cc}%
\frac{z}{1-z}r_{j}I_{2}+(\hat{w}-\hat{s})I_{2} & -\frac{1}{1-z}r_{j}%
\Delta_{2j}^{\prime}-\hat{u}I_{2}\\
-\frac{1}{1-z}r_{j}\Delta_{2j}-\hat{u}I_{2} & \frac{1}{1-z}I_{2}+\hat{v}I_{2}%
\end{array}
\right)  ,\label{Omegatilde}%
\end{align}
and%
\[
\delta_{j}=\left(  \hat{w}-\hat{s}\right)  \left(  1+\hat{v}-z\hat{v}\right)
+r_{j}\left(  \hat{u}+z\hat{v}-1\right)  -\left(  1-z\right)  \hat{u}^{2}.
\]

By definition,
\begin{align*}
\left\vert \hat{s}\right\vert  & \leq\frac{p}{T}\left\Vert D^{-1}\right\Vert
,\left\vert \hat{u}\right\vert \leq\frac{p}{T}\operatorname*{tr}\left\Vert
D^{-1}C^{\prime}M^{-1}\right\Vert ,\\
\left\vert \hat{v}\right\vert  & \leq\frac{p}{T}\left\Vert M^{-1}\right\Vert
,\text{ and }\left\vert \hat{w}\right\vert \leq\frac{p}{T}\operatorname*{tr}%
\left\Vert D^{-1}C^{\prime}M^{-1}CD^{-1}\right\Vert .
\end{align*}
In the proof of Lemma \ref{Rigour1}, we show that the norms $\left\Vert
D^{-1}\right\Vert ,$ $\left\Vert D^{-1}C^{\prime}\right\Vert ,$ and
$\left\Vert M^{-1}\right\Vert $ almost surely remain bounded as
$p,T\rightarrow_{c}\infty.$ Hence, $\hat{s},$ $\hat{u},$ $\hat{v},$ and
$\hat{w}$ are also almost surely bounded. Further, by definition,%
\[
r_{j}\Delta_{2j}=R_{j}-I_{2}\text{ and }r_{j}\Delta_{2j}^{\prime}%
=R_{j}^{\prime}-I_{2},
\]
where $R_{j}$ is an orthogonal matrix, so that $\left\Vert r_{j}\Delta
_{2j}\right\Vert $ and $\left\Vert r_{j}\Delta_{2j}^{\prime}\right\Vert $ are
clearly bounded uniformly in $j.$ Therefore, the norm of matrix $\tilde
{\Omega}_{j}$ almost surely remains bounded as $p,T\rightarrow_{c}\infty$,
uniformly in $j.$ Regarding $\delta_{j},$ which appear in the denominator on
the right hand side of (\ref{Omega}), the Supplementary Appendix establishes
the following result.

\begin{lemma}
\label{small delta lemma}There exists $\zeta>0$ such that, for any $z$ with
$\mathfrak{R}z=0$ and $\mathfrak{I}z>\zeta,$ almost surely,%
\[
\lim\inf_{p,T\rightarrow_{c}\infty}\max_{j=1,...,T/2}\left\vert \delta
_{j}\right\vert >c^{2}/\left(  1-c^{2}\right)  .
\]

\end{lemma}

The above results imply that, for $z$ with $\mathfrak{R}z=0$ and
$\mathfrak{I}z>\zeta,$ $\left\Vert \hat{\Omega}_{j}\right\Vert $ almost surely
remains bounded as $p,T\rightarrow_{c}\infty$, uniformly in $j.$ Therefore, by
Lemma \ref{Rigour1},%
\begin{equation}
\varepsilon_{(j)}^{\prime}M^{-1}\varepsilon_{(j)}=\hat{S}_{j}%
+o(1),\label{approximation1}%
\end{equation}
where $o(1)\overset{a.s.}{\rightarrow}0$ as $p,T\rightarrow_{c}\infty,$
uniformly in $j.$

A straightforward algebra reveals that%
\[
\hat{S}_{j}=\frac{\left(  \hat{w}-\hat{s}-r_{j}\right)  \hat{v}-\hat{u}^{2}%
}{\delta_{j}}.
\]
Using this in equations (\ref{approximation1}) and (\ref{germ of eq2}), we
obtain%
\begin{align*}
\hat{m}  & =\frac{2}{p}%
{\displaystyle\sum_{j=0}^{T/2-1}}
\frac{\left(  \hat{w}-\hat{s}-r_{j+1}\right)  \hat{v}-\hat{u}^{2}}%
{\delta_{j+1}}+o(1)\\
& =\frac{2}{p}%
{\displaystyle\sum_{j=1}^{T/2-1}}
\frac{f_{1}\left(  j\pi/T\right)  }{\left(  1-z\right)  f_{1}\left(
j\pi/T\right)  +f_{2}\left(  j\pi/T\right)  }+o(1),
\end{align*}
where, in the latter expression, the term corresponding to $j=0$ is included
in the $o(1)$ term to take into account the special definition of $r_{1}.$

As follows from Lemma \ref{small delta lemma} and the boundedness of $\hat
{s},\hat{u},\hat{v},$ and $\hat{w},$ the derivative
\[
\frac{\mathrm{d}}{\mathrm{d}\varphi}\frac{f_{1}\left(  \varphi\right)
}{\left(  1-z\right)  f_{1}\left(  \varphi\right)  +f_{2}\left(
\varphi\right)  }%
\]
almost surely remains bounded by absolute value as $p,T\rightarrow_{c}\infty,$
uniformly in $\varphi\in\left[  0,2\pi\right]  .$ Therefore%
\[
\frac{2}{p}%
{\displaystyle\sum_{j=1}^{T/2-1}}
\frac{f_{1}\left(  j\pi/T\right)  }{\left(  1-z\right)  f_{1}\left(
j\pi/T\right)  +f_{2}\left(  j\pi/T\right)  }=\frac{2}{\pi c}\int_{0}^{\pi
/2}\frac{f_{1}\left(  \varphi\right)  \mathrm{d}\varphi}{\left(  1-z\right)
f_{1}\left(  \varphi\right)  +f_{2}\left(  \varphi\right)  }+o(1).
\]
The statement of Proposition \ref{the second equation} now follows by noting
that the latter integral is one quarter of the integral over $\left[
0,2\pi\right]  .\square$

A similar analysis of equation (\ref{equation4}) gives us another proposition,
describing $\hat{m}$ as function of $\tilde{s},\tilde{u},\tilde{v},$ and
$\tilde{w},$ where%
\begin{align*}
\tilde{s}  & =\frac{1}{T}\operatorname*{tr}\left[  A^{-1}\right]  ,\text{
}\tilde{u}=\frac{1}{T}\operatorname*{tr}\left[  A^{-1}C\tilde{M}^{-1}\right]
,\\
\tilde{v}  & =\frac{1}{T}\operatorname*{tr}\left[  \tilde{M}^{-1}\right]
,\text{ and }\\
\tilde{w}  & =\frac{1}{T}\operatorname*{tr}\left[  A^{-1}C\tilde{M}%
^{-1}C^{\prime}A^{-1}\right]  .
\end{align*}
We omit the proof because it is very similar to that of Proposition
\ref{the second equation}.

\begin{proposition}
\label{the fourth equation}There exists $\zeta>0$ such that, for any $z$ with
$\mathfrak{R}z=0$ and $\mathfrak{I}z>\zeta,$ we have%
\begin{equation}
\hat{m}=\frac{1}{2\pi c}\int_{0}^{2\pi}\frac{g_{1}}{\left(  1-z\right)
g_{1}+g_{2}\left(  \varphi\right)  }\mathrm{d}\varphi+o(1),\text{
where}\label{asyeq2}%
\end{equation}%
\begin{align*}
g_{1}  & =\left(  \tilde{w}-\tilde{s}-1\right)  \tilde{v}-\tilde{u}^{2},\\
g_{2}\left(  \varphi\right)   & =\tilde{v}-4\sin^{2}\varphi\left(  \tilde
{s}+1-\tilde{u}-\tilde{w}\right)  ,
\end{align*}
and $o(1)\overset{a.s}{\rightarrow}0$, as $p,T\rightarrow_{c}\infty.$
\end{proposition}

Although we now have two asymptotic equations for $\hat{m},$ (\ref{asyeq1})
and (\ref{asyeq2}), they contain many unknowns: $\hat{s},\hat{u},\hat{v}%
,\hat{w},$ and the corresponding variables with tildes. The following result
establishes simple relationships between the unknowns with hats and tildes.

\begin{lemma}
\label{connections}We have the following three identities%
\begin{equation}
\hat{u}=\tilde{u},\text{ }z\tilde{v}+\hat{s}=\hat{w},\text{ and }z\hat
{v}+\tilde{s}=\tilde{w}.\label{connection identities}%
\end{equation}

\end{lemma}

\paragraph{Proof of Lemma \ref{connections}.}

The identity $\hat{u}=\tilde{u}$ is established by the following sequence of
equalities%
\begin{align*}
T\hat{u}  & =\operatorname*{tr}D^{-1}C^{\prime}M^{-1}=\operatorname*{tr}%
D^{-1}C^{\prime}\left(  CD^{-1}C^{\prime}-zA\right)  ^{-1}\\
& =\operatorname*{tr}\left(  C-zA\left(  C^{\prime}\right)  ^{-1}D\right)
^{-1}=\operatorname*{tr}\left(  C^{\prime}-zD\left(  C\right)  ^{-1}A\right)
^{-1}\\
& =\operatorname*{tr}A^{-1}C\left(  C^{\prime}A^{-1}C-zD\right)
^{-1}=\operatorname*{tr}A^{-1}C\tilde{M}^{-1}=T\tilde{u}.
\end{align*}
The relationship $z\tilde{v}+\hat{s}=\hat{w}$ is obtained as follows%
\begin{align*}
T\left(  z\tilde{v}+\hat{s}\right)   & =\operatorname*{tr}D^{-1}\left(
zI_{p}\left(  C^{\prime}A^{-1}CD^{-1}-zI_{p}\right)  ^{-1}+I_{p}\right) \\
& =\operatorname*{tr}D^{-1}\left(  -I_{p}+C^{\prime}A^{-1}CD^{-1}\left(
C^{\prime}A^{-1}CD^{-1}-zI_{p}\right)  ^{-1}+I_{p}\right) \\
& =\operatorname*{tr}D^{-1}\left(  I_{p}-DC^{-1}A\left(  C^{\prime}\right)
^{-1}z\right)  ^{-1}\\
& =\operatorname*{tr}D^{-1}C^{\prime}\left(  CD^{-1}C^{\prime}-Az\right)
^{-1}CD^{-1}=T\hat{w}.
\end{align*}
The identity $z\hat{v}+\tilde{s}=\tilde{w}$ is obtained similarly to
$z\tilde{v}+\hat{s}=\hat{w}$ by interchanging the roles of $D,C$ and
$A,C^{\prime}$.$\square$

The identities (\ref{connection identities}) imply the following equality%
\[
\left(  1-z\right)  f_{1}\left(  \varphi\right)  +f_{2}\left(  \varphi\right)
=\left(  1-z\right)  g_{1}+g_{2}\left(  \varphi\right)  .
\]
We denote the reciprocal of the common value of the right and left hand sides
of this equality as $\hat{h}\left(  z,\varphi\right)  .$ A direct calculation
shows that%
\begin{equation}
\hat{h}\left(  z,\varphi\right)  =\left(  \left(  1-z\right)  \left(
z\tilde{v}\hat{v}-\hat{u}^{2}\right)  +z\tilde{v}+4\sin^{2}\varphi\left(
z\hat{v}+\hat{u}-1\right)  \right)  ^{-1},\label{hfunction}%
\end{equation}
and the asymptotic relationships (\ref{asyeq1}) and (\ref{asyeq2}) can be
written in the following form%
\begin{equation}
\left\{
\begin{array}
[c]{l}%
\hat{m}=\frac{1}{2\pi c}\int_{0}^{2\pi}\hat{h}\left(  z,\varphi\right)
\left(  \left(  z\tilde{v}-4\sin^{2}\varphi\right)  \hat{v}-\hat{u}%
^{2}\right)  \mathrm{d}\varphi+o(1)\\
\hat{m}=\frac{1}{2\pi c}\int_{0}^{2\pi}\hat{h}\left(  z,\varphi\right)
\left(  \left(  z\hat{v}-1\right)  \tilde{v}-\hat{u}^{2}\right)
\mathrm{d}\varphi+o(1)
\end{array}
\right.  .\label{system of two}%
\end{equation}
This can be viewed as an asymptotic system of two equations with four
unknowns: $\hat{m},\tilde{v},\hat{v},$ and $\hat{u}.$ We shall now complete
the system by establishing the other two asymptotic relationships connecting
these unknowns.

Multiplying both sides of the identity%
\begin{equation}
MA^{-1}=CD^{-1}C^{\prime}A^{-1}-zI_{p}\label{another identity}%
\end{equation}
by $AM^{-1},$ taking trace, dividing by $p$, and rearranging terms, we obtain%
\begin{equation}
1+z\hat{m}=\frac{1}{p}\operatorname*{tr}\left[  CD^{-1}C^{\prime}%
M^{-1}\right]  .\label{equation1}%
\end{equation}
Next, we analyze (\ref{equation1}) similarly to the above analysis of
(\ref{equation2}). That is, first, we note that%
\begin{equation}
\frac{1}{p}\operatorname*{tr}\left[  CD^{-1}C^{\prime}M^{-1}\right]  =\frac
{1}{p}%
{\displaystyle\sum_{j=1}^{T/2}}
\operatorname*{tr}\left[  \Delta_{2j}^{\prime}\varepsilon_{(j)}^{\prime}%
D^{-1}C^{\prime}M^{-1}\varepsilon_{(j)}\right]  .\label{germ of eq1}%
\end{equation}
Then elementary algebra, based on the Sherman-Morrison-Woodbury formula
(\ref{identity}), yields%
\begin{align}
\varepsilon_{(j)}^{\prime}D^{-1}C^{\prime}M^{-1}\varepsilon_{(j)}  &
=r_{j}\left(  r_{j}I_{2}+s_{j}\right)  ^{-1}s_{j}\Delta_{2j}\left(  v_{j}-
\left[  v_{j},u_{j}^{\prime}\right]  \Omega_{j}\left[  v_{j},u_{j}^{\prime
}\right]  ^{\prime}\right) \label{stem of eq1}\\
& +r_{j}\left(  r_{j}I_{2}+s_{j}\right)  ^{-1}\left(  u_{j}-\left[
u_{j},w_{j}\right]  \Omega_{j}\left[  v_{j},u_{j}^{\prime}\right]  ^{\prime
}\right)  .\nonumber
\end{align}
Multiplying both sides of (\ref{stem of eq1}) by $\Delta_{2j}^{\prime}$ and
replacing $s_{j},u_{j},v_{j},$ and $w_{j}$ by $\hat{s}I_{2},\hat{u}I_{2}%
,\hat{v}I_{2},$ and $\hat{w}I_{2},$ respectively, yields an asymptotic
approximation to $\Delta_{2j}^{\prime}\varepsilon_{(j)}^{\prime}%
D^{-1}C^{\prime}M^{-1}\varepsilon_{(j)},$ which can be used in
(\ref{germ of eq1}) and (\ref{equation1}) to produce the following result. Its
proof, as well as the proof of (\ref{stem of eq1}), are given in the
Supplementary Appendix.

\begin{proposition}
\label{the first equation}There exists $\zeta>0$ such that, for any $z$ with
$\mathfrak{R}z=0$ and $\mathfrak{I}z>\zeta,$ we have%
\begin{equation}
1+z\hat{m}=\frac{1}{2\pi c}\int_{0}^{2\pi}\hat{h}\left(  z,\varphi\right)
\left(  2\hat{u}\sin^{2}\varphi+z\tilde{v}\hat{v}-\hat{u}^{2}\right)
\mathrm{d}\varphi+o(1),\text{ where}\label{asyeq3}%
\end{equation}
$o(1)\overset{a.s}{\rightarrow}0$, as $p,T\rightarrow_{c}\infty.$
\end{proposition}

One might think that the remaining asymptotic relationship can be obtained by
using the identity
\begin{equation}
\tilde{M}D^{-1}=C^{\prime}A^{-1}CD^{-1}-zI_{p},\label{parallel identity}%
\end{equation}
which parallels (\ref{another identity}). Unfortunately, following this idea
delivers a relationship equivalent to (\ref{asyeq3}). Therefore, instead of
using (\ref{parallel identity}), we consider the identity%
\begin{equation}
\frac{1}{p}\operatorname*{tr}\left[  C^{\prime}M^{-1}\right]  =\frac{1}%
{p}\operatorname*{tr}\left[  DD^{-1}C^{\prime}M^{-1}\right]
,\label{last identity}%
\end{equation}
which yields%
\begin{equation}
\frac{1}{p}%
{\displaystyle\sum_{j=1}^{T/2}}
\operatorname*{tr}\left[  \Delta_{2j}\varepsilon_{(j)}^{\prime}M^{-1}%
\varepsilon_{(j)}\right]  =\frac{1}{p}%
{\displaystyle\sum_{j=1}^{T/2}}
\operatorname*{tr}\left[  \Delta_{1j}\varepsilon_{(j)}^{\prime}D^{-1}%
C^{\prime}M^{-1}\varepsilon_{(j)}\right]  .\label{germ of eq4}%
\end{equation}
Then, we proceed as in the above analysis of (\ref{germ of eq1}) and
(\ref{germ of eq2}) to obtain the remaining asymptotic relationship. The proof
of the following proposition is given in the Supplementary Appendix.

\begin{proposition}
\label{the last equation}There exists $\zeta>0$ such that, for any $z$ with
$\mathfrak{R}z=0$ and $\mathfrak{I}z>\zeta,$ we have%
\begin{equation}
0=\frac{1}{2\pi c}\int_{0}^{2\pi}\hat{h}\left(  z,\varphi\right)  \left(
4\hat{v}\sin^{2}\varphi+2\hat{u}\right)  \mathrm{d}\varphi+o(1),\text{
where}\label{asyeq4}%
\end{equation}
$o(1)\overset{a.s}{\rightarrow}0$, as $p,T\rightarrow_{c}\infty.$
\end{proposition}

Summing up the results in Propositions \ref{the second equation},
\ref{the fourth equation}, \ref{the first equation}, and
\ref{the last equation}, the unknowns $\hat{m},\hat{v},\tilde{v},$ and
$\hat{u}$ must satisfy the following system of asymptotic equations%
\begin{equation}
\left\{
\begin{array}
[c]{l}%
\hat{m}=\frac{1}{2\pi c}\int_{0}^{2\pi}\hat{h}\left(  z,\varphi\right)
\left(  \left(  z\tilde{v}-4\sin^{2}\varphi\right)  \hat{v}-\hat{u}%
^{2}\right)  \mathrm{d}\varphi+o(1)\\
\hat{m}=\frac{1}{2\pi c}\int_{0}^{2\pi}\hat{h}\left(  z,\varphi\right)
\left(  \left(  z\hat{v}-1\right)  \tilde{v}-\hat{u}^{2}\right)
\mathrm{d}\varphi+o(1)\\
1+z\hat{m}=\frac{1}{2\pi c}\int_{0}^{2\pi}\hat{h}\left(  z,\varphi\right)
\left(  2\hat{u}\sin^{2}\varphi+z\tilde{v}\hat{v}-\hat{u}^{2}\right)
\mathrm{d}\varphi+o(1)\\
0=\frac{1}{2\pi c}\int_{0}^{2\pi}\hat{h}\left(  z,\varphi\right)  \left(
4\hat{v}\sin^{2}\varphi+2\hat{u}\right)  \mathrm{d}\varphi+o(1)
\end{array}
\right.  .\label{Asymptotic system}%
\end{equation}

\subsubsection{Solving the system}

Recall that the unknowns $\hat{m},\hat{v},\tilde{v},$ and $\hat{u}$ in the
asymptotic relationships (\ref{Asymptotic system}) depend on $p,T.$ The
definition (\ref{mhat def}) of $\hat{m}$ implies that $\left\vert \hat
{m}\right\vert $ is bounded by $\left(  \mathfrak{I}z\right)  ^{-1}.$ Further,
as shown in the proof of Proposition \ref{the second equation}, $\hat{u}$ and
$\hat{v}$ are a.s. bounded by absolute value, and it can be similarly shown
that $\tilde{v}$ is a.s. bounded by absolute value. Therefore, there exist a
subsequence of $p,T$ along which $\hat{m},\hat{v},\tilde{v},$ and $\hat{u}$
a.s. converge to some limits $m,v,y,$ and $u.$

These limits must satisfy a non-asymptotic system of equations%
\begin{equation}
\left\{
\begin{array}
[c]{l}%
m=\frac{1}{2\pi c}\int_{0}^{2\pi}h\left(  z,\varphi\right)  \left(  \left(
zy-4\sin^{2}\varphi\right)  v-u^{2}\right)  \mathrm{d}\varphi\\
m=\frac{1}{2\pi c}\int_{0}^{2\pi}h\left(  z,\varphi\right)  \left(  \left(
zv-1\right)  y-u^{2}\right)  \mathrm{d}\varphi\\
1+zm=\frac{1}{2\pi c}\int_{0}^{2\pi}h\left(  z,\varphi\right)  \left(
2u\sin^{2}\varphi+zvy-u^{2}\right)  \mathrm{d}\varphi\\
0=\frac{1}{2\pi c}\int_{0}^{2\pi}h\left(  z,\varphi\right)  \left(  2v\sin
^{2}\varphi+u\right)  \mathrm{d}\varphi
\end{array}
\right.  ,\label{the system}%
\end{equation}
where%
\[
h\left(  z,\varphi\right)  =\left[  \left(  1-z\right)  \left(  zvy-u^{2}%
\right)  +zy+4\sin^{2}\varphi\left(  zv+u-1\right)  \right]  ^{-1}.
\]
Let us consider, until further notice, only such $z$ that have zero real part,
$\mathfrak{R}z=0$, and the imaginary part satisfying $\mathfrak{I}z>\zeta,$
for some $\zeta>0$. Let us solve system (\ref{the system}) for $m.$ Adding two
times the last equation to the first one, and subtracting the second equation
we obtain%
\begin{equation}
0=\frac{1}{2\pi c}\int_{0}^{2\pi}h\left(  z,\varphi\right)  \left(
y+2u\right)  \mathrm{d}\varphi.\label{zero}%
\end{equation}

Note that $\int_{0}^{2\pi}h\left(  z,\varphi\right)  \mathrm{d}\varphi\neq0.$
Otherwise, from the second equation of (\ref{the system}), we have $m=0,$
which cannot be true because $\hat{m}$ is the Stieltjes transform of the
empirical distribution of the squared canonical correlations, all of which lie
between zero and one. Indeed, clearly, for any $0\leq\lambda\leq1$ and $z$
with $\mathfrak{R}z=0,$%
\[
\mathfrak{I}\left(  \frac{1}{\lambda-z}\right)  =\frac{\mathfrak{I}z}%
{\lambda^{2}+\left(  \mathfrak{I}z\right)  ^{2}}\geq\frac{\mathfrak{I}%
z}{1+\left(  \mathfrak{I}z\right)  ^{2}}.
\]
Therefore, $\mathfrak{I}\hat{m}\geq\mathfrak{I}z/\left(  1+\left(
\mathfrak{I}z\right)  ^{2}\right)  ,$ and $\hat{m}$ cannot converge to $m=0$.

Since $\int_{0}^{2\pi}h\left(  z,\varphi\right)  \mathrm{d}\varphi\neq0,$
(\ref{zero}) yields
\begin{equation}
y+2u=0\label{reduction1a}%
\end{equation}
with $y\neq0$ and $u\neq0$ (if one of them equals zero, the other equals zero
too, and $m=0$ by the second equation of (\ref{the system}), which is
impossible). Since $u\neq0,$ the last equation implies that $v\neq0$ as well.

Further, subtracting from the third equation the sum of $z$ times the second
and $u/v$ times the last equation, and using (\ref{reduction1a}), we obtain%
\begin{equation}
1=\frac{1}{2\pi c}\int_{0}^{2\pi}h\left(  z,\varphi\right)  \frac{u}{v}\left(
2zv+u\right)  \left(  zv-v-1\right)  \mathrm{d}\varphi.\label{germ of red2}%
\end{equation}
This equation, together with the second equation of (\ref{the system}) yield%
\begin{equation}
m=\frac{v\left(  2zv+u-2\right)  }{\left(  1+v-zv\right)  \left(
2zv+u\right)  }.\label{reduction2}%
\end{equation}

Next, for the integrand in the last equation of (\ref{the system}), we have%
\begin{align}
& h\left(  z,\varphi\right)  \left(  2v\sin^{2}\varphi+u\right)  =\frac{1}%
{2}\frac{v}{zv+u-1}\label{hsimplified}\\
& +h\left(  z,\varphi\right)  \frac{u}{2}\left(  \frac{\left(  1-z\right)
v\left(  2zv+u\right)  +2\left(  2zv+u-1\right)  }{zv+u-1}\right)  .\nonumber
\end{align}
This assumes that
\begin{equation}
zv+u-1\neq0.\label{non-degeneracy}%
\end{equation}
If not, then
\[
h\left(  z,\varphi\right)  =\left[  \left(  1-z\right)  \left(  zvy-u^{2}%
\right)  +zy\right]  ^{-1}%
\]
would not depend on $\varphi$ and the last equation of (\ref{the system})
would imply that $u+v=0.$ The latter equation and the equality $zv+u-1=0$
would yield $v=-\left(  1-z\right)  ^{-1},$ which when combined with the
second equation of (\ref{the system}) would give us $m=-c^{-1}\left(
1-z\right)  ^{-1},$ which cannot be true because $m,$ being a limit of
$\hat{m}$, must satisfy $\mathfrak{I}m\geq0$ for $\mathfrak{I}z>0.$

Equations (\ref{germ of red2}), (\ref{hsimplified}), and the last equation of
(\ref{the system}) imply that%
\begin{equation}
u=\frac{2c}{2c-1-\left(  1-z\right)  v\left(  1-c\right)  }%
-2zv.\label{reduction3}%
\end{equation}
Combining this with (\ref{reduction2}) yields%
\begin{equation}
m=v\frac{1-c}{c}.\label{reduction4}%
\end{equation}

Finally, elementary calculations given in the Supplementary Appendix show that%
\begin{equation}
\left(  \frac{1}{2\pi}\int_{0}^{2\pi}\frac{1}{x+2\sin^{2}\varphi}%
\mathrm{d}\varphi\right)  ^{2}=\frac{1}{x\left(  x+2\right)  },\label{fact}%
\end{equation}
where $x\in\mathbb{C}\backslash\left[  -2,0\right]  $. Using (\ref{fact}),
(\ref{germ of red2}), and the definition of $h\left(  z,\varphi\right)  $, we
obtain the following relationship%
\begin{align}
& \left(  \frac{2cv\left(  zv+u-1\right)  }{u\left(  2zv+u\right)  \left(
zv-v-1\right)  }\right)  ^{2}\label{ugly}\\
& =\frac{4\left(  zv+u-1\right)  ^{2}}{u\left(  \left(  1-z\right)  \left(
-2zv-u\right)  -2z\right)  \left(  -u+uz+2\right)  \left(  u+2vz-2\right)
},\nonumber
\end{align}
that holds as long as%
\[
\frac{u\left(  \left(  1-z\right)  \left(  -2zv-u\right)  -2z\right)
}{2\left(  zv+u-1\right)  }\in\mathbb{C}\backslash\left[  -2,0\right]  .
\]
The latter inclusion holds because otherwise $h\left(  z,\varphi\right)  $ is
not a bounded function of $\varphi,$ which would contradict Lemma
\ref{small delta lemma}.

Using (\ref{reduction3}) in (\ref{ugly}), and simplifying, we find that there
exist only three possibilities. Either%
\begin{equation}
v=-\frac{1}{1-z},\label{possibility1}%
\end{equation}
or%
\begin{equation}
1-\left(  c+cz-1\right)  v+z\left(  1-z\right)  \left(  1-c\right)
v^{2}=0,\label{possibility2}%
\end{equation}
or%
\begin{equation}
\frac{c}{1-c}-\left(  c+cz-z\right)  v+z\left(  1-z\right)  \left(
1-c\right)  v^{2}=0.\label{possibility3}%
\end{equation}

Equation (\ref{possibility1}) cannot hold because otherwise, (\ref{reduction4}%
) would imply that $\mathfrak{I}m<0,$ which is impossible as argued above.
Equation (\ref{possibility2}) taken together with (\ref{reduction3}) implies
that%
\[
u+zv-1=0,
\]
which was ruled out above. This leaves us with (\ref{possibility3}), so that,
using (\ref{reduction4}), we get%
\begin{equation}
m=\frac{-\left(  z-c-cz\right)  \pm\sqrt{\left(  z-c-cz\right)  ^{2}-4c\left(
1-z\right)  z}}{2z\left(  1-z\right)  c}.\label{second possibility}%
\end{equation}

For $z\in\mathbb{C}^{+}$ with $\mathfrak{R}z=0,$ the imaginary part of the
right hand side of (\ref{second possibility}) is negative when `$-$' is used
in front of the square root. Here we choose the branch of the square root,
with the cut along the positive real semi-axis, which has positive imaginary
part. Since $\mathfrak{I}m$ cannot be negative, we conclude that%
\begin{equation}
m=\frac{-\left(  z-c-cz\right)  +\sqrt{\left(  z-c-cz\right)  ^{2}-4c\left(
1-z\right)  z}}{2z\left(  1-z\right)  c}.\label{the limit}%
\end{equation}
But the right hand side of the above equality is the value of the limit of the
Stieltjes transforms of the eigenvalues of the multivariate beta matrix
$B_{p}\left(  p,\left(  T-p\right)  /2\right)  $ as $p,T\rightarrow_{c}%
\infty.$ This can be verified directly by using the formula for such a limit,
given for example in Theorem 1.6 of Bai, Hu, Pan and Zhou (2015). As follows
from Wachter (1980), the weak limit of the empirical distribution of the
eigenvalues of the multivariate beta matrix $B_{p}\left(  p,\left(
T-p\right)  /2\right)  $ as $p,T\rightarrow_{c}\infty$ equals $W\left(
\lambda;c/\left(  1+c\right)  ,2c/\left(  1+c\right)  \right)  $.

Equation (\ref{the limit}) shows that, for $z$ with $\mathfrak{R}z=0$ and
$\mathfrak{I}z>\zeta,$ any converging subsequence of $\hat{m}$ converges to
the same limit. Hence, $\hat{m}$ a.s. converges for all $z$ with
$\mathfrak{R}z=0$ and $\mathfrak{I}z>\zeta.$ Note that $\hat{m}$ is a sequence
of bounded analytic functions in the domain $\left\{  z:\mathfrak{I}%
z>\delta\right\}  ,$ where $\delta$ is an arbitrary positive number.
Therefore, by Vitaly's convergence theorem (see Titchmarsh (1939), p.168)
$\hat{m}$ a.s. converges to $m,$ described by (\ref{the limit}), for any
$z\in\mathbb{C}^{+}.$ The almost sure convergence of $\hat{F}_{p,T}\left(
\lambda\right)  $ (and thus, also of $F_{p,T}\left(  \lambda\right)  $) to the
Wachter distribution follows from the Continuity Theorem for the Stieltjes
transforms (see, for example, Corollary 1 in Geronimo and Hill (2003)).

\subsection{Proof of Theorem \ref{Levy}}

First, let us show that the weak limit $F_{0}\left(  \lambda\right)  $ of
$F_{\gamma}\left(  \lambda\right)  $ as $\gamma\rightarrow0$ exists and equals
the continuous part of the Marchenko-Pastur distribution with density
(\ref{MP density}). By definition and Theorem \ref{main}, $F_{\gamma}\left(
\lambda\right)  $ is the (scaled) Wachter d.f. $W\left(  \gamma\lambda
;\gamma/\left(  1+\gamma\right)  ,2\gamma/\left(  1+\gamma\right)  \right)  .$
Therefore, by (\ref{densityW}) and (\ref{support boundaries}), the density,
$f_{\gamma}(\lambda)$, and the boundaries of the support, $\left[  \hat{b}%
_{-},\hat{b}_{+}\right]  ,$ of the distribution $F_{\gamma}$ equal
\[
f_{\gamma}(\lambda)=\frac{1+\gamma}{2\pi}\frac{\sqrt{\left(  \hat{b}%
_{+}-\lambda\right)  \left(  \lambda-\hat{b}_{-}\right)  }}{\lambda\left(
1-\gamma\lambda\right)  },\text{ and}%
\]%
\[
\hat{b}_{\pm}=\left(  \sqrt{2}\mp\sqrt{1-\gamma}\right)  ^{-2}.
\]

As $\gamma\rightarrow0,$ $\hat{b}_{\pm}\rightarrow a_{\pm},$ where $a_{\pm
}=\left(  1\pm\sqrt{2}\right)  ^{2}$ as in (\ref{MP boundaries}), and
$f_{\gamma}(\lambda)$ converges to the density given by (\ref{MP density}).
This implies the weak convergence of $F_{\gamma}\left(  \lambda\right)  $ to
$F_{0}\left(  \lambda\right)  $ with $F_{0}$ supported on $\left[  a_{-}%
,a_{+}\right]  $ and having density (\ref{MP density}).

To establish the theorem, it remains to show that, as $p\rightarrow\infty$,
$F_{p,\infty}(\lambda)$ weakly converges to $F_{0}(\lambda),$ in probability.
Recall that the weak convergence is metrized by the L\'{e}vy distance
$\mathcal{L}\left(  \cdot,\cdot\right)  $. We need to show that for any
$\delta>0,$ there exists $p_{0}$ such that (s.t.) for all $p>p_{0},$%
\begin{equation}
\Pr\left(  \mathcal{L}\left(  F_{0},F_{p,\infty}\right)  <\delta\right)
>1-\delta.\label{need to show}%
\end{equation}

Let $\gamma>0$ be so small that%
\begin{equation}
\mathcal{L}\left(  F_{0},F_{\gamma}\right)  <\delta/4.\label{Levy1}%
\end{equation}
For any $p,$ let $T_{\gamma}$ be the smallest even integer satisfying
$p/T_{\gamma}\leq\gamma.$ That is,%
\[
T_{\gamma}=\min_{T\in2\mathbb{Z}}\left\{  T:p/T\leq\gamma\right\}  .
\]
For any $T_{\infty}>T_{\gamma},$ by the triangle inequality, we have%
\begin{equation}
\mathcal{L}\left(  F_{0},F_{p,\infty}\right)  \leq\mathcal{L}\left(
F_{0},F_{\gamma}\right)  +\mathcal{L}\left(  F_{\gamma},F_{p,T_{\gamma}%
}\right)  +\mathcal{L}\left(  F_{p,T_{\gamma}},F_{p,T_{\infty}}\right)
+\mathcal{L}\left(  F_{p,T_{\infty}},F_{p,\infty}\right)  ,\label{triangle}%
\end{equation}
where $F_{p,T_{\gamma}}$ and $F_{p,T_{\infty}}$ denote the empirical
distributions of eigenvalues of%
\begin{equation}
\frac{T}{p}CD^{-1}C^{\prime}A^{-1},\label{RawMatrixAgain}%
\end{equation}
with $T=T_{\gamma}$ and $T=T_{\infty},$ respectively.

By Theorem \ref{main}, $\mathcal{L}\left(  F_{\gamma},F_{p,T_{\gamma}}\right)
$ a.s. converges to zero as $p\rightarrow\infty.$ Therefore, for all
sufficiently large $p$, we have%
\begin{equation}
\Pr\left(  \mathcal{L}\left(  F_{\gamma},F_{p,T_{\gamma}}\right)
<\delta/4\right)  >1-\delta/4.\label{Levy2}%
\end{equation}
Further, as shown by Johansen (1988, 1991), for any $p,$ as $T_{\infty
}\rightarrow\infty,$ the eigenvalues of (\ref{RawMatrixAgain}) with
$T=T_{\infty}$ jointly converge in distribution to those of%
\begin{equation}
\frac{1}{p}\int_{0}^{1}\left(  \mathrm{d}B\right)  B^{\prime}\left(  \int%
_{0}^{1}BB^{\prime}\mathrm{d}u\right)  ^{-1}\int_{0}^{1}B\left(
\mathrm{d}B\right)  ^{\prime}.\label{JohansenLimitAgain}%
\end{equation}
Therefore, for any $p$ and all sufficiently large $T_{\infty},$ we have%
\begin{equation}
\Pr\left(  \mathcal{L}\left(  F_{p,T_{\infty}},F_{p,\infty}\right)
<\delta/4\right)  >1-\delta/4.\label{Levy4}%
\end{equation}

Let us denote the sum of $\mathcal{L}\left(  F_{0},F_{\gamma}\right)  ,$
$\mathcal{L}\left(  F_{\gamma},F_{p,T_{\gamma}}\right)  ,$ and $\mathcal{L}%
\left(  F_{p,T_{\infty}},F_{p,\infty}\right)  $ as $\mathcal{L}_{\gamma
,p,T_{\infty}}.$ By (\ref{triangle}), we have%
\begin{equation}
\mathcal{L}\left(  F_{0},F_{p,\infty}\right)  \leq\mathcal{L}_{\gamma
,p,T_{\infty}}+\mathcal{L}\left(  F_{p,T_{\gamma}},F_{p,T_{\infty}}\right)
.\label{two terms}%
\end{equation}
Inequalities (\ref{Levy1}), (\ref{Levy2}), and (\ref{Levy4}) show that for any
$\delta>0,$ there exists $\gamma_{\delta}>0$ such that (s.t.) for any positive
$\gamma<\gamma_{\delta},$ there is a $p_{\gamma}$ s.t. for any $p>p_{\gamma},$
there is a $T_{p}$ s.t. for any $T_{\infty}>T_{p}$%
\begin{equation}
\Pr\left(  \mathcal{L}_{\gamma,p,T_{\infty}}<3\delta/4\right)  >1-\delta
/2.\label{Levy124}%
\end{equation}
The subscripts in $\gamma_{\delta},$ $p_{\gamma}$ and $T_{p}$ signify
dependence on the value of the corresponding parameter. Inequalities
(\ref{Levy124}) and (\ref{two terms}) would establish (\ref{need to show}) as
long as we are able to show that for any $\delta>0,$ there exists
$\tilde{\gamma}_{\delta}>0$ s.t. for any positive $\gamma<\tilde{\gamma
}_{\delta}, $ there is a $\tilde{p}_{\gamma}$ s.t. for any $p>\tilde
{p}_{\gamma}$ \textit{and any} $\tilde{T}_{p},$ there exists $T_{\infty
}>\tilde{T}_{p}$ s.t.%
\begin{equation}
\Pr\left(  \mathcal{L}\left(  F_{p,T_{\gamma}},F_{p,T_{\infty}}\right)
<\delta/4\right)  >1-\delta/2.\label{Levy3}%
\end{equation}

Let us denote $\xi=\sqrt{T}\varepsilon,$ where $\varepsilon$ is a $p\times T$
matrix with i.i.d. $N(0,1/T)$ entries, as defined in Section 2. We shall
assume that, as $p,T$ change, $\xi$ represents $p\times T$ sections of a fixed
infinite array of i.i.d. standard normal random variables. Consider%
\[
M_{p,T}=\frac{T}{p}\left(  \frac{\xi\xi^{\prime}}{T}\right)  ^{-1/2}\frac
{\xi\Delta_{2}^{\prime}\xi^{\prime}}{T}\left(  \frac{\xi\Delta_{1}\xi^{\prime
}}{T}\right)  ^{-1}\frac{\xi\Delta_{2}\xi^{\prime}}{T}\left(  \frac{\xi
\xi^{\prime}}{T}\right)  ^{-1/2}.
\]
So defined matrix $M_{p,T}$ is identical to the real symmetric matrix
$\frac{T}{p}A^{-1/2}CD^{-1}C^{\prime}A^{-1/2}.$ The above definition is
formulated in terms of $\xi$ to clarify that $M_{p,T}$ depends on $T$ not only
via the term $T/p,$ but also through $A,C,$ and $D.$ Note that $F_{p,T_{\gamma
}}$ and $F_{p,T_{\infty}}$ are the empirical distributions of eigenvalues of
$M_{p,T_{\gamma}}$ and $M_{p,T_{\infty}},$ respectively. The following lemma
is established in the Supplementary Appendix.

\begin{lemma}
\label{alpha12} For any $\tau>0$ there exists $\gamma_{\tau}>0$ s.t. for any
positive $\gamma<\gamma_{\tau}$, there is a $\tilde{p}_{\gamma}$ s.t. for any
$p>\tilde{p}_{\gamma}$ and any $\tilde{T}_{p},$ there exists $T_{\infty
}>\tilde{T}_{p}$ s.t. with probability larger than $1-\tau,$ $M_{p,T_{\gamma}%
}-M_{p,T_{\infty}}$ can be represented as the sum of two real symmetric
matrices $S$ and $R$,%
\[
M_{p,T_{\gamma}}-M_{p,T_{\infty}}=S+R,
\]
where $\left\Vert S\right\Vert \leq K\sqrt{\gamma},$ $\operatorname*{rank}%
R\leq\tau p$, and $K$ is an absolute constant.
\end{lemma}

Finally, let $F_{SR}$ be the empirical distribution of eigenvalues of
$M_{p,T_{\gamma}}-S=M_{p,T_{\infty}}+R.$ Then, by Theorem A45 (norm
inequality) of Bai and Silverstein (2010),%
\[
\mathcal{L}\left(  F_{p,T_{\gamma}},F_{SR}\right)  \leq\left\Vert S\right\Vert
\leq K\sqrt{\gamma},
\]
whereas by their Theorem A43 (rank inequality),%
\[
\mathcal{L}\left(  F_{SR},F_{p,T_{\infty}}\right)  \leq\frac{1}{p}%
\operatorname*{rank}R\leq\tau.
\]
Therefore, by Lemma \ref{alpha12} and the triangle inequality, for any
$\tau>0$ there exists $\gamma_{\tau}>0$ s.t. for any positive $\gamma
<\gamma_{\tau}$, there is a $\tilde{p}_{\gamma}$ s.t. for any $p>\tilde
{p}_{\gamma}$ and any $\tilde{T}_{p},$ there exists $T_{\infty}>\tilde{T}_{p}$
s.t.%
\[
\Pr\left(  \mathcal{L}\left(  F_{p,T_{\gamma}},F_{p,T_{\infty}}\right)
<\tau+K\sqrt{\gamma}\right)  >1-\tau.
\]
For $\tau=\delta/8,$ this inequality implies (\ref{Levy3}) with $\tilde
{\gamma}_{\delta}=\min\left\{  \gamma_{\tau},\left(  \delta/8K\right)
^{2}\right\}  .$ Combining (\ref{Levy3}) with (\ref{Levy124}) yields
(\ref{need to show}), which completes the proof.

\end{document}